\documentclass[11pt]{amsart}

\usepackage{amsmath,amsthm, amscd, amssymb, amsfonts}
\usepackage{epsfig}
\input xy
\xyoption{all}

\newcommand{\Hc}{{\mathcal H}}

\newcommand{\otb}{{\overline{\otimes}}}
\newcommand{\otk}{{\otimes_{\ku}}}
\newcommand{\Mo}{{\mathcal M}}

\newcommand{\nic}{{\mathfrak B}}

\newcommand{\Ss}{{\mathcal S}}
\newcommand{\ot}{{\otimes}}

\newcommand{\Ac}{{\mathcal A}}

\newcommand{\ca}{{\mathcal C}}

\newcommand{\Do}{{\mathcal D}}
\newcommand{\Bc}{{\mathcal B}}

\newcommand{\Fc}{{\mathcal F}}
\newcommand{\Gc}{{\mathcal G}}
\newcommand{\YD}{{\mathcal YD}}

\newcommand{\ku}{{\Bbbk}}

\newcommand{\Na}{{\mathbb N}}

\newcommand{\uno}{{\mathbf 1}}

\newcommand{\id}{\mbox{\rm id\,}}
\newcommand{\Tot}{\mbox{\rm Tot\,}}

\newcommand{\Der}{\mbox{\rm Der\,}}

\newcommand\Rep{\operatorname{Rep}}

\newcommand\co{\operatorname{co}}

\newcommand\Hom{\operatorname{Hom}}

\newcommand\Homk{\operatorname{Hom}_{\ku}}
\newcommand{\End}{\operatorname{End}}

\newcommand{\gr}{\mbox{\rm gr\,}}

\newcommand{\grl}{\mbox{\rm gr\,}_{\! l}}

\newcommand{\grc}{\mbox{\rm gr\,}_{\! c}}

\renewcommand{\_}[1]{\mbox{$_{\left( #1 \right)}$}}

\theoremstyle{plain}

\numberwithin{equation}{section}

\newtheorem{teo}{Theorem}[section]

\newtheorem{lema}[teo]{Lemma}

\newtheorem{cor}[teo]{Corollary}

\newtheorem{prop}[teo]{Proposition}

\newtheorem{claim}{Claim}[section]

\theoremstyle{definition}

\newtheorem{defi}[teo]{Definition}

  \newtheorem{exa}[teo]{Example}

\theoremstyle{remark}

\newtheorem{rmk}[teo]{Remark}

\def\pf{\begin{proof}}

\def\epf{\end{proof}}

\theoremstyle{remark}

\begin{document}

\title[Module categories over  pointed
Hopf algebras] {Module categories over pointed Hopf algebras}
\author[Mombelli]{
Mart\'\i n Mombelli }
\thanks{This work was conducted during a post-doc fellowship at the Ludwig
Maximilians Universit\"{a}t, M\"{u}nchen, granted by the Deutscher
Akademischer Austausch Dienst (DAAD)}
\address{Mathematisches Institut, Universit\"{a}t M\"{u}nchen, Theresienstra\ss e 39,
\newline \indent D-80333 M\"{u}nchen, Germany.}
\email{mombelli@mate.uncor.edu, martin10090@gmail.com
\newline \indent\emph{URL:}\/ http://www.mate.uncor.edu/mombelli}
\begin{abstract} We develop some techniques for studying  exact
module categories over some families of pointed finite-dimensional
Hopf algebras. As an application we classify exact module
categories over the tensor category of representations of the
small quantum groups $u_q(\mathfrak{sl}_2)$.
\end{abstract}


\date{\today}
\maketitle

\setcounter{tocdepth}{1}

\section{Introduction}

Given a tensor category $\ca$, a module category over $\ca$ is an
abelian category $\Mo$ together with an exact functor $\ca\times
\Mo\to \Mo$ satisfying natural associativity and unity axioms.
This rather general notion appears, and is implicitly present, in
many different areas of mathematics and mathematical physics, such
as the theory of weak Hopf algebras \cite{O1}, subfactor theory
\cite{BEK}, extensions of vertex algebras \cite{KO}, affine Hecke
algebras \cite{BO} and conformal field theory, see for example
\cite{FS}, \cite{CS1}, \cite{CS2}. In this case one is interested
in the tensor category of representations of a certain vertex
algebra, and irreducible objects of a module category are
interpreted as boundary conditions for the conformal field theory.

\medbreak

The language of module categories was used with profit in diverse
papers. It is present the study of fusion categories, see
\cite{ENO1}, \cite{ENO2}, and in relation with dynamical twists
\cite{M}. In \cite{N} module categories were used to describe some
properties of semisimple Hopf algebras.

\medbreak

Etingof and Ostrik proposed \cite{eo}, as an important problem,
the cla\-ssification of a certain class of module categories,
called \emph{exact},  and not only of semisimple module
categories. This was done in the framework of the study of finite
tensor categories not necessarily semisimple, motivated by its
relation, among others, with logarithmic conformal field theories.
Such classification is available only for few examples, see
\cite{eo}, \cite{O1}, \cite{O2}, \cite{O3}.

\medbreak

This paper arises mainly  motivated by the problem of classifying
module categories over the tensor category of representations of
the small quantum groups $u_q(\mathfrak{g})$, for a simple Lie
algebra $\mathfrak{g}$  and $q$ a primitive root of unity. One of
the goals of this paper is the classification in the case
$\mathfrak{g}=\mathfrak{sl}_2$. This particular problem was
proposed by V. Ostrik at the meeting \emph{Groupes quantiques
dynamiques et categories de fusion}, at CIRM, Marseille 2008.

\medbreak

The paper is organized as follows. In section \ref{preliminaries}
we give a brief account of the basic definitions of graded and
coradically graded Hopf algebras, twisting in Hopf algebras and
Hopf-Galois extensions. In section \ref{mod-cate} we recall some
basic facts on exact module categories as introduced by Etingof
and Ostrik \cite{eo}, and  we also recall some results on module
categories over the tensor category of representations of a
finite-dimensional Hopf algebra. In particular we recall the main
result from \cite{AM} that any exact indecomposable module
category over $\Rep(H)$, $H$ a finite-dimensional Hopf algebra, is
equivalent to the category of left $K$-modules, for some right
$H$-simple left $H$-comodule $K$.

\medbreak

Section \ref{section:graded-alg} is devoted to the study of
filtered (graded) comodule algebras over a filtered (graded) Hopf
algebra. In section \ref{section:lift} we study the relation
between \emph{liftings} of graded comodule algebras, that is
comodule algebras $A$ such that $\gr A\simeq G$, where $G$ is a
given graded comodule algebra, and formal deformations of $G$.
This is completely analogous to the study of liftings of Hopf
algebras in \cite{DCY}.

The different liftings of a given comodule algebra are controlled
by a certain double complex which we present in section
\ref{section:bicomplex} generalizing previous results of D. Yau
\cite{Y2}, \cite{y1}.

\medbreak

In section \ref{section:comod-pointed} we show that if $A$ is a
certain graded comodule algebra over a Hopf algebra $H$
constructed by a bosonization as $H=\nic(V)\# \ku G$, where
$\nic(V)$ is the Nichols algebra of a Yetter-Drinfeld module $V$
over the group algebra $\ku G$ over a cyclic group $G$, then $A$
is determined by a certain homogeneous coideal subalgebra in
$\nic(V)$ and a subgroup of $G$.

\medbreak

Using the techniques developed  previously, in section
\ref{section:classification} we present the classification of
exact module categories over the Taft Hopf algebras, the Radford
Hopf algebras, the book Hopf algebras and over the Lusztig's small
quantum groups $u_q(\mathfrak{sl}_2)$.

\section{Preliminaries and notation }\label{preliminaries}

Throughout this work $\ku$ will denote an algebraically closed
field of characteristic zero. All vector spaces, algebras and
linear categories will be considered over $\ku$. By $\ku^{\times}$
we shall denote the non-zero elements of $\ku$.

\medbreak

If $A$ is an algebra, we shall denote by ${}_A\Mo$ the category of
finite-dimensional left $A$-modules. If $H$ is a Hopf algebra
$\Rep(H)$ will denote the tensor category of finite-dimensional
representations of $H$ and ${}^H_H\YD$ will denote the category of
Yetter-Drinfeld modules of $H$.

\medbreak

Given a Hopf algebra $H$, a left $H$-comodule algebra is a
collection $(A, \mu,\lambda)$, where $(A,\mu)$ is an algebra,
$\lambda:A\to H\otk A$ is the $H$-comodule structure on $A$ such
that $\lambda$ is an algebra map. We shall sometimes denote this
comodule algebra by $(A,\lambda)$. We will denote it by $(A,
\mu,\lambda)$ only when special emphasis on the product is needed.
The set of coinvariants of $A$ is $A^{\co H}=\{a\in A:
\lambda(a)=1\ot a\}.$
 \medbreak

If $A$ is an $H$-comodule algebra via $\lambda:A\to H\otk A$, we
shall say that a (right) ideal $J$ is $H$-costable if
$\lambda(J)\subseteq H\otk A$. We shall say that $A$ is (right)
$H$-simple if there is no nontrivial $H$-costable (right) ideal in
$A$.

\medbreak

Let $A$ be a left $H$-comodule algebra. We denote by
${}^{H}\Mo_{A}$ the category of left $H$-comodules, right
$A$-modules $M$ such that the right $A$-module structure $ M\otk
A\to M$ is an $H$-comodule map. If $P\in{}^{H}\Mo_{A}$ we will
denote by $\End^H_A(P)$ the space of $A$-module and $H$-comodule
endomorphisms.

\medbreak

A left coideal subalgebra of $H$ is a subalgebra $K\subseteq H$
such that $\Delta(K)\subseteq H\otk K$. It is well-known that if
$H$ is finite-dimensional, $K$ is a left coideal subalgebra of $H$
and $Q$ denotes the coalgebra $H/ HK^+$, then  the pair of
functors $\mathcal{F}:{}^{H}\Mo_{K}\to {}^Q\Mo$,
$\mathcal{F}(M)=M/ MK^+$ and $\mathcal{G}: {}^Q\Mo\to
{}^{H}\Mo_{K}$, $\mathcal{G}(V)=H\square_Q V$ give an equivalence
of categories. See \cite{Ma}, \cite{Sk}.

\medbreak
 For a  coalgebra $C$ we shall denote by $C_0$ the
coradical and by $C_0\subseteq C_1 \subseteq C_1 \dots \subseteq
C_n\subseteq \dots$ the coradical filtration. For more details see
\cite{Sw}.

\medbreak

Let us recall that a coalgebra $C$ is said to be
$\Na_0$-\emph{graded}, or \emph{graded} for short, if $C=
\bigoplus_{i} C(i)$, such that $\Delta(C(n))\subseteqq
\bigoplus^n_{i=0} C(i)\otk C(n-i)$ for all $n\in \Na_0$ and
$\varepsilon(C(n))=0$ for all $n\geq 1$. A coalgebra is
\emph{coradically graded} if it is graded and the coradical
filtration is given by $C_n=\bigoplus^n_{i=0} C(i)$.

A Hopf algebra $H$ is $\Na_0$-\emph{graded}, or \emph{graded} for
short, if it is $\Na_0$-graded as an algebra and a coalgebra by
the same grading.

\medbreak

A \emph{ filtration} on a Hopf algebra $H$ is an algebra
filtration $H^0\subseteq H^1 \subseteq \dots \subseteq H^m=H$ such
that for all $n=0\dots m$
$$ \Delta(H^n)\subseteq \sum^n_{i=0} \; H^i \otk H^{n-i}.$$

If $H$ is a filtered Hopf algebra such that $H^0$ is a Hopf
subalgebra then the graded algebra associated to the filtration $
\gr H= \bigoplus_{n\geq 0} H^n/ H^{n-1}$, here $ H^{-1}=0$, is a
 graded Hopf algebra. The associated graded Hopf algebra to the coradical
 filtration will be denoted by $\grc H$. It is a well-known fact that
$\grc H$ is a coradically graded Hopf algebra.

 \medbreak

Let $H$ be a graded Hopf algebra. A \emph{lifting} of $H$ is
filtered Hopf algebra structure $U$ on the same underlying vector
space $H$ such that $\gr U=H$.

Two liftings $U$ and $V$ of $H$ are said to be \textit{equivalent}
if there is a Hopf algebra isomorphism $\phi:U\to V$ such that
$\gr \phi=\id_H$. The set of equivalence classes of liftings of
$H$ will be denoted by $Lift(H)$.

\subsection{Twisting comodule algebras}

Let us recall that a Hopf 2-cocycle is a  map $\sigma: H\otk H\to
\ku$, invertible with respect to convolution,  such that
\begin{align}\label{2-cocycle}
\sigma(x\_1, y\_1)\sigma(x\_2y\_2, z) &= \sigma(y\_1,
z\_1)\sigma(x, y\_2z\_2),
\\
\label{2-cocycle-unitario} \sigma(x, 1) &= \varepsilon(x) =
\sigma(1, x),
\end{align}
for all $x,y, z\in H$. Using this cocycle there is a new Hopf
algebra structure constructed over the same coalgebra $H$ with the
product described by
$$
x._{[\sigma]}y = \sigma(x_{(1)}, y_{(1)}) \sigma^{-1}(x_{(3)},
y_{(3)})\, \, x_{(2)}y_{(2)}, \qquad x,y\in H.
$$
We shall denote by  $H^{\sigma}$ this new Hopf algebra. We shall
denote by ${}_\sigma H$ the space $H$ with new multiplication
given by
\begin{align}\label{twisted-left} x\cdot y=\sigma(x\_2, y\_2)\; x\_1 y\_1,
\end{align}
for all $x,y\in H$. The algebra ${}_\sigma H$ is a left
$H$-comodule algebra with coaction given by the coproduct of $H$.

\medbreak

We shall say that $\sigma$ is \emph{cocentral} if for all $x, y\in
H$
\begin{equation}\label{2-cocycle-central}
\sigma(x_{(1)}, y_{(1)})\, x_{(2)}y_{(2)} = \sigma(x_{(2)},
y_{(2)})\, x_{(1)}y_{(1)}.\end{equation}

It is obvious that $\sigma$ is cocentral if and only if
$H^{\sigma}=H$. A useful way to check if a certain 2-cocycle
$\sigma$ is cocentral is to prove that the coproduct ${}_\sigma
H:\to {}_\sigma H\otk H$ makes ${}_\sigma H$ a right $H$-comodule
algebra. This observation will be used later without further
mention.

\medbreak

If $\sigma: H\otimes H\to \ku$ is a 2-cocycle and $K$ is a left
$H$-comodule algebra, then we can define a new product in $K$ by
\begin{align}\label{sigma-product} a._{\sigma}b = \sigma(a_{(-1)},
b_{(-1)})\, a_{(0)}.b_{(0)},
\end{align}
$a,b\in K$. We shall denote by $K_{\sigma}$ this new algebra. It
is easy to see that $K_{\sigma}$ is a left $H^{\sigma}$-comodule
algebra. The following basic result will be useful later.

\begin{lema}\label{equivalencia-cocentral} Let $\sigma: H\otimes H\to \ku$
be a Hopf 2-cocycle and $K$ be a left $H$-comodule algebra. There
is an equivalence of categories ${}^{H}\Mo_{K}\simeq
{}^{H^{\sigma}}\Mo_{K_{\sigma}}$.

\end{lema}
\pf If $M\in {}^{H}\Mo_{K}$ we define a right $K_{\sigma}$-module
structure by
$$  m\cdot_{\sigma} k= \sigma(m\_{-1}, k\_{-1})\, m\_0\cdot
k\_0,$$ for all $k\in K$, $m\in M$. We shall denote by
$M_{\sigma}$ the object $M$ with this new action and the same left
$H$-comodule structure. By a straightforward computation one can
prove that $M_{\sigma}\in {}^{H^{\sigma}}\Mo_{K_{\sigma}}$ and
that the functor $M\longmapsto M_{\sigma}$ is an equivalence of
categories. \epf In particular if $K\subseteq H$ is a left coideal
subalgebra, $Q=H/HK^+$ and $\sigma: H\otimes H\to \ku$ is a
cocentral 2-cocycle then the categories
${}^{H\!}\Mo_{K_{\sigma}}$, ${}^Q\Mo$ are equivalent.

\medbreak

The following definition will be used later.

\begin{defi} A Hopf 2-cocycle $\sigma: H\otimes H\to \ku$ is \emph{compatible}
with a left coideal subalgebra $K\subseteq H$ if for any $x,y\in
K$, $\sigma(x\_2, y\_2) \; x\_1y\_1\in K$.
\end{defi}
In the case when $\sigma$ is compatible with $K$ we denote
${}_\sigma K$ the left $H$-comodule algebra obtained from $K$
using multiplication given as in \eqref{twisted-left}. Any
cocentral 2-cocycle is compatible with any left coideal subalgebra
and in this case ${}_\sigma K=K_\sigma.$

\subsection{Hopf Galois extensions} Let $H$ be a Hopf algebra. A
left $H$-Galois extension of an algebra $R$ is an algebra $A$ with
a left $H$-comodule algebra structure, $A\to H\otk A$, $a\to
a\_{-1}\ot a\_0$, such that $A^{\co H}=R$ and the canonical map
$can: A\ot_R A \to H\ot A$, $can(a\ot b)=a\_{-1}\ot\, a\_0b$ is a
bijection. For more details on this subject the reader is referred
to \cite{SS}.

\medbreak

It is well known that if $\sigma:H\otk H\to \ku$ is a 2-cocycle
then ${}_{\sigma}H$ is a $H$-Galois extension of the field $\ku$.
In the finite-dimensional case the converse also holds, that is,
if $H$ is finite-dimensional and $A$ is a finite-dimensional left
$H$-Galois extension of the field then there is a 2-cocycle
$\sigma$ such that $A\simeq {}_{\sigma}H$.

\section{Exact module categories over finite tensor
categories}\label{mod-cate}

We recall the definitions of module categories over finite tensor
categories. For more details the reader is referred to \cite{O1},
\cite{eo}.

\medbreak

A {\em module category} over a tensor category $\ca$ is an abelian
category $\Mo$ equipped with an exact bifunctor $\otimes: \ca
\times \Mo \to \Mo$ and natural associativity and unit
isomorphisms $m_{X,Y,M}: (X\otimes Y)\otimes M \to X\otimes
(Y\otimes M)$, $\ell_M: \uno \otimes M\to M$ such that for any
$X,Y,Z\in \ca$, $M\in \Mo$.
\begin{align}\label{action1}
(\id\otimes m_{Y,Z,M})m_{X,Y\otimes Z,M}(a_{X,Y,Z}\otimes id) &=
m_{X,Y,Z\otimes M}\, m_{X\otimes Y,Z,M},
\\\label{action2}
(\id\otimes \ell_M)m_{X,\uno ,Y} &=r_X\otimes \id.
\end{align}

  A {\em module functor} between module
categories $\Mo$ and $\Mo'$ over a tensor category $\ca$ is a pair
$(\Fc,c)$, where $\Fc:\Mo \to \Mo'$ is a $\ku$-linear functor and
$c_{X,M}: \Fc(X\otimes M)\to X\otimes \Fc(M)$ is a natural
isomorphism such that  for any $X, Y\in \ca$, $M\in \Mo$:
\begin{align}\label{modfunctor1}
(\id_X\otimes c_{Y,M})c_{X,Y\otimes M}\Fc(m_{X,Y,M}) &=
m_{X,Y,\Fc(M)}\, c_{X\otimes Y,M},
\\\label{modfunctor2}
\ell_{\Fc(M)} \,c_{\uno ,M} &=\Fc(\ell_{M}).
\end{align}

Let $\Mo_1$ and $\Mo_2$ be module categories over $\ca$. We shall
denote by $\Hom_{\ca}(\Mo_1, \Mo_2)$ the category whose objects
are module functors $(\Fc, c)$ from $\Mo_1$ to $\Mo_2$. A morphism
between  $(\Fc,c)$ and $(\Gc,d)\in\Hom_{\ca}(\Mo_1, \Mo_2)$ is a
natural transformation $\alpha: \Fc \to \Gc$ such that for any
$X\in \ca$, $M\in \Mo_1$:
\begin{gather}
\label{modfunctor3} d_{X,M}\alpha_{X\otimes M} =
(\id_{X}\ot\alpha_{M})c_{X,M}.
\end{gather}

 Two module categories $\Mo_1$ and $\Mo_2$ over
$\ca$ are {\em equivalent} if there exist module functors
$F:\Mo_1\to \Mo_2$ and $G:\Mo_2\to \Mo_1$ and natural isomorphisms
$\id_{\Mo_1} \to F\circ G$, $\id_{\Mo_2} \to G\circ F$ that
satisfy \eqref{modfunctor3}.

 The {\em direct sum} of two module categories
$\Mo_1$ and $\Mo_2$ over a tensor category $\ca$  is the
$\ku$-linear category $\Mo_1\times \Mo_2$ with coordinate-wise
module structure.

 A module category is {\em indecomposable} if it
is not equivalent to a direct sum of two non trivial module
categories.

\medbreak In this paper we further assume that all module
categories have finitely many isomorphism classes of simple
objects. The following definition of a very important class of
module categories is due to P. Etingof and V. Ostrik.

\begin{defi}\label{exactmc} \cite{eo}
A module category $\Mo$ over a finite tensor category $\ca$ is
\emph{exact} if  for any projective $P\in \ca$ and any $M\in \Mo$,
the object $P\ot M$ is projective in $\Mo$.
\end{defi}

Any semisimple finite module category over a finite tensor
category $\ca$ is exact. A direct sum of finite module categories
is exact if and only if each summand is exact. Therefore, any
exact module category over $\ca$ is a finite direct product of
exact indecomposable module categories over $\ca$, see \cite[Prop.
3.9]{eo}.

\subsection{Module categories over Hopf algebras}

We recall some  results concerning exact module categories over
the category of representations of  Hopf algebras. Let $H$ be a
finite-dimensional Hopf algebra.

\medbreak

If $\lambda:A\to H\otk A$ is a left $H$-comodule algebra then the
category of finite-dimensional left $A$-modules ${}_A\Mo$ is a
module category over $\Rep(H)$ with action $\otb:\Rep(H)\times
{}_A\Mo\to {}_A\Mo$, $X\otb M=X\otk M$, for all $X\in \Rep(H),
M\in {}_A\Mo$. The left $A$-module structure on $X\otk M$ is given
by $\lambda$, that is, if $a\in A$, $x\in X, m\in M$ then $a\cdot
(x\ot m)= a\_{-1}\cdot x\ot\; a\_{0}\cdot m$.

\medbreak

If $A$ is right $H$-simple then ${}_A\Mo$ is an exact
indecomposable module category over $\Rep(H)$, see \cite[Prop.
1.20]{AM}. The exactness is a consequence of the beautiful results
obtained by Skryabin \cite{Sk} on comodule algebras.

\begin{rmk} In particular the category ${}_A\Mo$ is an exact indecomposable
module category over $\Rep(H)$ when
\begin{enumerate}
    \item $A=K_{\sigma}$, $K\subseteq H$ is a coideal subalgebra
    and $\sigma: H\otk H\to \ku$ is a cocentral
Hopf 2-cocycle,
    \item or, more generally, $A= {}_\sigma K$
    where $K\subseteq H$ is a left coideal subalgebra and $\sigma:H\otk H
    \to \ku$ is a 2-cocycle compatible with $K$.
\end{enumerate}

\end{rmk}

\begin{teo}\cite[Theorem 3.3]{AM}\label{mod-overhopf} If $\Mo$ is an exact idecomposable
module category over $\Rep(H)$ then $\Mo\simeq {}_A\Mo$ for some
right $H$-simple left comodule algebra $A$ with $A^{\co
H}=\ku$.\qed
\end{teo}

Two left $H$-comodule algebras $A$ and $B$ are \emph{Morita
equivariant equivalent}, and we shall denote it by $A\sim_M B$, if
the module categories ${}_A\Mo$, ${}_B\Mo$ are equivalent as
module categories over $\Rep(H)$.

\begin{prop}\label{eq1}\cite[Prop. 1.24]{AM} The algebras $A$ and $B$ are
Morita equivariant equivalent if and only if there exists $P\in
{}^{H\!}\Mo_{B}$ such that $A\simeq \End_B(P)$ as $H$-comodule
algebras.\qed
\end{prop}

The left $H$-comodule structure on $\End_B(P)$ is given by
$\lambda:\End_B(P)\to H\otk \End_B(P)$, $\lambda(T)=T\_{-1}\ot
T\_0$ where
\begin{equation}\label{h-comod} \langle\alpha, T\_{-1}\rangle\,
T_0(p)=\langle\alpha, T(p\_0)\_{-1}\Ss^{-1}(p\_{-1})\rangle\,
T(p\_0)\_0,\end{equation} for any $\alpha\in H^*$,
$T\in\End_B(P)$, $p\in P$. It is easy to prove that
$\End_B(P)^{\co H}= \End^H_B(P)$.

\medbreak

Let $K\subseteq H$ be a left coideal subalgebra and $\sigma:H\otk
H\to \ku$ be a cocentral 2-cocycle. Let us denote by $Q$ the
quotient $H/ HK^+$.
\begin{lema}\label{morita-equivariantness} Let $A$ be a right $H$-simple left $H$-comodule
algebra with trivial coinvariants Morita equivariant equivalent to
$K_{\sigma}$. Then there exists an indecomposable object $P\in
{}^Q\Mo$ such that $A\simeq \End_{K_{\sigma}}((H\square_Q
P)_{\sigma})$ and $\dim A= (\dim P)^2 \dim K$.\end{lema}

\pf The existence of the object $P$ such that $A\simeq
\End_{K_{\sigma}}((H\square_Q P)_{\sigma})$ follows by Proposition
\ref{eq1} and Lemma \ref{equivalencia-cocentral}. Since the
coinvariants are trivial $A^{\co H}=\ku$ then $
\End^H_B(H\square_Q P)=\ku$ and this implies that $H\square_Q P$
must be indecomposable thus $P$ is also indecomposable. Since
$H\simeq K\ot Q$ as left $K$-modules, right $Q$-comodules, then
$H\square_Q P\simeq K\ot P$, whence $\dim A= (\dim P)^2 \dim
K$.\epf

Let us further assume that $Q$ is a pointed cosemisimple
coalgebra. Thus, if $P\in {}^Q\Mo$ is an indecomposable object
then is 1-dimensional and there exists a group-like element $g\in
Q$ and $0\neq v\in P$ such that $P=\ku<v>$ and $\delta(v)=g\ot v$
is the left $Q$-comodule structure on $P$. In this case
$H\square_Q P\simeq gK$.

\begin{lema}\label{morita-cocentral} Under the above assumptions
if $A$ is a left $H$-comodule algebra such that $A\sim_M
K_{\sigma}$ then $A\simeq (gKg^{-1})_{\sigma^g}$, for some
group-like element $g\in G(H)$ and $\sigma^g:H\otk H\to \ku$ is
the Hopf 2-cocycle given by
$$\sigma^g(x, y)=\sigma(g^{-1}xg,g^{-1}yg), $$
for any $x,y\in H$.
\end{lema}
\pf From the above considerations $A\simeq
\End_{K_{\sigma}}((gK)_{\sigma})$. Any element $T\in
\End_{K_{\sigma}}((gK)_{\sigma})$ is determined by the value
$T(g)$. Indeed, any element $gx\in gK$ can be written uniquely
$gx=g\cdot_{\sigma} \big(\sigma^{-1}(g,x\_1)\; x\_2\big)$. Here,
for any $x,y\in H$ we are denoting $x\cdot_{\sigma}
y=\sigma(x\_1,y\_1)\, x\_2y\_2$ the product as in
\eqref{sigma-product}.

The maps $\phi:
 (gKg^{-1})_{\sigma^g}\to \End_{K_{\sigma}}((Kg)_{\sigma})$,
 $\psi: \End_{K_{\sigma}}((Kg)_{\sigma})\to  (gKg^{-1})_{\sigma^g}$
given by
$$\phi(gy^{-1}g)(g\cdot x)=g\cdot (x\cdot y), \quad \psi(T)=g^{-1}T(g),$$ for any $x, y\in K$ $T\in
 \End_{K_{\sigma}}((Kg)_{\sigma})$ are well-defined maps one the inverse of
 each other. It is straightforward to verify that both maps are
 $H$-comodule morphisms.

\epf

\begin{exa} Let $G$ be a finite group, $F\subseteq G$ be a subgroup
and $H=\ku G$ is the group algebra. The group algebra $\ku F$ is a
left coideal subalgebra of $H$. Let also $\psi\in
Z^2(G,\ku^{\times})$ be a (necessarily cocentral) 2-cocycle. Let
us denote by $\ku^{\psi} F$ the twisted group algebra. An
indecomposable object $P\in {}^Q\Mo$, $Q=\ku G/ (\ku F)^+\ku G$
must be one-dimensional, say $P=<\ku\overline{g}>$. If $A\sim_M
\ku^{\psi} F$ then there exists an element $g\in G$ such that
$A\simeq \ku^{\psi^g} F^g$, where $F^g$ is the conjugate of $F$.
This result is contained in \cite[Thm 2]{O1}.
\end{exa}

\section{Graded and filtered comodule
algebras}\label{section:graded-alg}

Let us assume that $H$ is a finite-dimensional filtered Hopf
algebra with filtration given by $H^0\subseteq H^1 \subseteq\dots
\subseteq H^m=H$ such that $H^0$ is a Hopf subalgebra.

\medbreak

Given  a left $H$-comodule algebra $(A, \lambda)$, we shall say
that a filtration $A^0\subseteq A^1\subseteq\dots \subseteq A^m$
is a \emph{comodule algebra filtration} compatible with the
filtration of $H$ if it is an algebra filtration and for each
$n=0,\dots,m$
\begin{equation}\label{filtr-comod} \lambda(A^n)\subseteq \sum^n_{i=0}\, H^i \otk A^{n-i}.
\end{equation} Observe
that $A^0$ is a left $H^0$-comodule algebra. Whenever no confusion
arises we shall simply say that $A$ is a \emph{filtered}
$H$-comodule algebra.

\medbreak

Any filtration on the Hopf algebra induces a compatible filtration
on a comodule algebra.

\begin{lema}\label{induced-filt} Let $H^0\subseteq H^1 \subseteq\dots \subseteq H^m=H$
be a filtration on the Hopf algebra $H$ and $(A, \lambda)$ a
comodule algebra. Define $A^n=\lambda^{-1}(H^n\ot A)$ for all
$n=0\dots m,$ then $A^0\subseteq\dots \subseteq A^m=A$ is a
filtration compatible with the filtration on $H$.
\end{lema}
\pf Since $\lambda$ is an algebra map, for any $i, j$ we have that
$A^i A^j\subseteq A^{i+j}$. Let be $0\leq n\leq m$ and $a\in A^n$,
then $\lambda(a)\in H^n\ot A$, hence $$(\id_H\ot\lambda)
\lambda(a)= (\Delta\ot \id_A)\lambda(a)\in \sum^n_{i=0}\, H^i\ot
H^{n-i}\ot A.$$ Whence $\lambda(a)\in\sum^n_{i=0}\, H^i\ot
A^{n-i}$.

\epf As a consequence any comodule algebra $(A, \lambda)$ has a
distinguished filtration compatible with the coradical filtration
on $H$, \emph{the Loewy series} on $A$; it is the filtration
$A_0\subseteq A_1\subseteq\dots \subseteq A_m$ defined by
$A_n=\lambda^{-1}(H_n\ot A)$ for all $n\in \Na$.

\begin{defi}  Let $H=\oplus^m_{i=0} H(i)$ be a  graded Hopf algebra. We
shall say that a left $H$-comodule algebra $G$, with comodule
structure given by $\lambda:G\to H\otk G$, graded as an algebra
$G=\oplus^m_{i=0} G(i)$ is a \emph{graded comodule algebra} if for
each $0\leq n\leq m$
$$ \lambda(G(n))\subseteq \bigoplus^m_{i=0} H(i)\otk G(n-i).$$
A graded comodule algebra $G=\oplus^m_{i=0} G(i)$ is
\emph{Loewy-graded} if it is a graded comodule algebra and the
Loewy series is given by $G_n=\oplus^n_{i=0} G(i)$.

\smallbreak

Let $(A, \lambda)$ be a filtered left $H$-comodule algebra,
$A^0\subseteq A^1\subseteq\dots \subseteq A^m$ compatible with the
filtration $H^0\subseteq H^1 \subseteq\dots \subseteq H^m=H$. We
can consider the graded algebra $\gr A=\bigoplus_{n\geq 0} A^n/
A^{n-1}$, where $A^{-1}=0$. The graded algebra associated to the
Loewy series will be denoted by $\grl A$.

\smallbreak

There is a well defined map $\overline{\lambda}:\gr A\to \gr
H\ot\, \gr A$ such that the following diagram commutes
$$
\begin{CD}
 A^n@>{\lambda}>> \bigoplus^n_{i=0} H^i\otk A^{n-i}
\\
@VVV @VVV \\
 A^n/ A^{n-1} @>>> \big(\sum^n_{i=0} H^i\otk
 A^{n-i}\big)/ \sum^{n-1}_{j=0} H^j\otk
 A^{n-1-j}
 \\
@VVV @VV\simeq V \\
\gr A(n) @>{\overline{\lambda}}>> \sum^n_{i=0}\, \gr H(i) \otk\,
\gr A(n-i).
\end{CD}
$$

\begin{lema} The space $(\gr A, \overline{\lambda})$ is a graded
$ \gr H$-comodule algebra and $(\grl A, \overline{\lambda})$ is a
Loewy-graded $\grc H$-comodule algebra.
\end{lema}
\pf This is a well-known fact. The proof can be found for example
in \cite{AD}, \cite{Sw}.\epf

The following result will be useful when we relate module
categories over a Hopf algebra $H$ and over $\grc H$.
\begin{prop}\label{equivalence-grad}
Let $H^0\subseteq H^1 \subseteq\dots \subseteq H^m=H$ be a
filtration on the Hopf algebra $H$ such that $H^0$ is a semisimple
Hopf subalgebra. Let  $A$ be a left $H$-comodule algebra and
consider the filtration $A^i=\lambda^{-1}(H^i\otk A)$ for any
$i=0,\dots,m$. The following assertions are equivalent.
\begin{enumerate}
    \item $A^0$ is right $H^0$-simple.
    \item $A$ is right $H$-simple.
\end{enumerate}
\end{prop}
\pf Let us assume that $A^0$ is right $H^0$-simple. Let
$J\subseteq A$ be a right ideal $H$-costable. Consider the
filtration on $J$, $J^0\subseteq J^1\subseteq \dots J^m=J$, given
by $J^i= \lambda^{-1}(H^i\otk J)$. The space $J^0\subseteq A^0$ is
a right ideal $H^0$-costable, thus $J^0=0$ or $J^0=A^0$. In the
second case $1\in J$ and thus $J=A$. Let us assume that $J^0=0$.
Follows from Lemma \ref{induced-filt} that $\lambda(J^n)\subseteq
\sum^n_{i=0}\, H^i \otk J^{n-i}$. Then $\lambda(J^1)\subseteq
H_0\otk J^1$ thus $J^1=J^0$. Arguing inductively follows that
$J^n=J^{n-1}$ for all $n$, thus $J=0$.

\medbreak

Assume $A$ is right $H$-simple. Since $H^0$ is semisimple, using
\cite[Thm. 3.1]{L}, we obtain that the Jacobson radical $J(A^0)$
is a $H^0$-costable ideal of $A^0$, hence $J(A^0)A$ is a right
ideal $H$-costable in $A$. Since the radical is nilpotent the
equality $J(A^0)A=A$ is impossible, therefore $J(A^0)A=0$ and
$A^0$ is semisimple. Thus, the inclusion $A^0\subseteq A$ splits,
that is, there is a left $A^0$-module $B\subseteq A$ such that
$A=A^0\oplus B$.  Now, let $J\subseteq A^0$ be a nonzero right
$H^0$-costable ideal. Then $JA= J\oplus JB\subseteq J\oplus B$.
Since $A$ is right $H$-simple we obtain that  $JA=A$ which implies
that $J=A^0$. \epf

Under the same hypothesis as Proposition \ref{equivalence-grad} we
have the following result.
\begin{cor}  The
following assertions are equivalent.
\begin{enumerate}
    \item $A^0$ is right $H^0$-simple.
    \item $\gr A$ is
right $\gr H$-simple. \item $A$ is right $H$-simple.
\end{enumerate}
\end{cor}

\pf (1) is equivalent to (3) is Proposition
\ref{equivalence-grad}. But the same proof works for $\gr A$ since
$(\gr A)^0=A^0$.\epf

\section{Liftings and formal deformation of comodule
algebras}\label{section:lift} For the rest of this section $(H,
\Delta, m)$ will denote a graded Hopf algebra with grading
$H=\oplus^m_{i=0} H(i)$.

\medbreak

Let $G$ be a graded $H$-comodule algebra. If $U$ is a lifting for
$H$, we shall say that $A$ is a \emph{lifting} for $G$ along $U$
if $A=G$ as vector spaces and $A$ is a filtered left $U$-comodule
algebra and $\gr A\simeq G$ as left $H$-comodule algebras.
\end{defi}

Two liftings $A, B$ of $G$ along $U$ are said to be equivalent if
there is an isomorphism of filtered $U$-comodule algebras
$\psi:A\to B$ such that $\gr \psi=\id_G$. We shall denote by
$Lift(A,U)$ the set of equivalence classes of liftings of $A$
along $U$.

\subsection{Formal bialgebra deformation}
Let us recall graded bialgebra deformation from \cite{DCY}. See
also \cite{MW}. For any $l\in \Na \cup \{\infty\}$ consider the
free $\ku[t]/ (t^{l+1})$-module $H[t]/ (t^{l+1})$ as a graded
vector space over $\ku$ by $\deg t=1$ and $\deg h=i$ if $h\in
H(i)$. If $l=\infty$, $H[t]/ (t^{l+1})=H[t]$ and $\ku[t]/
(t^{l+1})=\ku[t]$.

\medbreak

An $l$-deformation of $H$ is a collection $\big(H[t]/ (t^{l+1}),
m^l, \Delta^l\big)$, where
$$m^l:(H\otk H)[t]/ (t^{l+1})\to H[t]/ (t^{l+1}), \;\;
\Delta^l:H[t]/ (t^{l+1})\to (H\otk H)[t]/ (t^{l+1})$$ are $\ku[t]/
(t^{l+1})$-linear maps homogeneous of degree 0, such that $H[t]/
(t^{l+1})$ is a bialgebra. This implies that there are maps
$\Delta^l_s:H\to H\otk H$, $m^l_s:H\otk H\to H$ of degree $-s$
such that for all $x,y \in H$
$$ \Delta^l(x)= \sum^l_{s=0} \, \Delta^l_s(x)\, t^s, \quad m^l(x\ot y)=
 \sum^l_{s=0} \, m^l_s(x\ot y)\,t^s.$$
It is also required that $\Delta^l_0=\Delta$ and $m^l_0=m$.

\medbreak

Two $l$-deformations $\big(H[t]/ (t^{l+1}), m^l, \Delta^l\big)$
and $\big(H[t]/ (t^{l+1}), \widetilde{m}^l,
\widetilde{\Delta}^l\big)$ are isomorphic if there exists a
bialgebra isomorphism $\phi:H[t]/ (t^{l+1})\to H[t]/ (t^{l+1})$
homogeneous of degree 0. In this case there are maps $\phi_s:H\to
H$ homogeneous of degree $-s$ such that for any $h\in H$,
$\phi(h)= \sum^l_{s=0} \,\phi_s(h) t^s$. Here it is required that
$\phi_0=\id_H$.

The set of equivalence of $l$-deformations is denoted by
$Iso^l(H)$. When $l=\infty$, this set is simply denoted by
$Iso(H)$.

\medbreak
 The following result is \cite[Thm. 2.2]{DCY}.

\begin{teo}\label{main-dyc} There exists a natural bijection between $Lift(H)$ and
$Iso(H)$.
\end{teo}

\subsection{Formal deformation for comodule algebras}

Let us fix a $l$-deformation $\big(H[t]/ (t^{l+1}), m^l,
\Delta^l\big)$   of $H$ that we shall denote by $H_l$ for short.
Let also $(G,\mu,\lambda)$ be a graded left $H$-comodule algebra.

\begin{defi} An $l$-\emph{deformation} of $G$ is a collection
$\big(G[t]/ (t^{l+1}), \mu^l, \lambda^l\big)$ where $\mu^l:(G\otk
G)[t]/ (t^{l+1})\to G[t]/ (t^{l+1})$, $\lambda^l: G[t]/
(t^{l+1})\to (H\otk G)[t]/ (t^{l+1})$ are $\ku[t]/
(t^{l+1})$-linear maps homogeneous of degree zero such that $G[t]/
(t^{l+1})$ is a left $H_l$-comodule algebra. In particular there
are maps $\lambda^l_s:G\to H\otk G$, $\mu^l_s:G\otk G\to G$
homogeneous of degree $-s$ such that for any $a, b\in G$
\begin{align*} \lambda^l(a)=\sum^l_{s=0}\;\lambda^l_s(a) \,
t^s,\quad\; \mu^l(a\ot b)=\sum^l_{s=0}\;\mu^l_s(a\ot b)\, t^s.
\end{align*}
It is required that $\lambda^l_0=\lambda$ and $\mu^l_0=\mu$.
\end{defi}

\medbreak

The associativity of the product, the coassociativity and the
compatibility of the coaction of $G[t]/ (t^{l+1})$ with the
product imply that for any $0\leq n\leq l$
\begin{equation}\label{compatibilid-def1} \sum_{i+j=n}\; \mu^l_i (\mu^l_j\ot\id_G)= \sum_{i+j=n}\; \mu^l_i
(\id_G\ot \mu^l_j),
\end{equation}
\begin{equation}\label{compatibilid-def2} \sum_{i+j=n}\; (\id_H\ot
\lambda^l_i)\lambda^l_j=\sum_{i+j=n}\;
(\Delta^l_i\ot\id_G)\lambda^l_j,
\end{equation}
\begin{equation}\label{compatibilid-def3} \sum_{i+j=n}\; \lambda^l_i
\mu^l_j=\sum_{i+j+k+r=n}\; (m^l_i\ot\mu^l_j)\,\tau_{23}\,
(\lambda^l_k\ot\lambda^l_r).
\end{equation}

 Two $l$-deformations $\big(G[t]/ (t^{l+1}), \mu^l,
\lambda^l\big)$, $\big(\widetilde{G}[t]/ (t^{l+1}),
\widetilde{\mu}^l, \widetilde{\lambda}^l\big)$ are equivalent if
there exists an homogeneous map of degree zero $\phi:G[t]/
(t^{l+1})\to \widetilde{G}[t]/ (t^{l+1})$ such that $\phi$ is a
$\ku[t]/ (t^{l+1})$-linear isomorphism of $H_l$-comodule algebras.
In particular there are homogeneous maps $\phi_s:G\to
\widetilde{G}$ of degree $-s$ such that
$$\phi(a)=\sum^l_{s=0}\; \phi_s(a)\, t^s,$$
for all $a\in G$. We required that $\phi_0=\id_G$. Denote by
$Iso^l(G,H_l)$ the set of equivalence classes of $l$-deformations
of $G$. The proof of the following Lemma is completely analogous
to \cite[Lemma 2.1]{DCY}.

\begin{lema}\label{restrictions} For $l>l'$ there exist a
restriction map $res_{l,l'}:Iso^l(G,H_l)\to Iso^{l'}(G,H_{l'})$,
and maps $r_l:Iso^{\infty}(G,H_{\infty})\to Iso^l(G,H_l)$ such
that
$$Iso^{\infty}(G,H_{\infty})\simeq
\underleftarrow{\rm lim}_{l\in \Na} \; Iso^l(G,H_l)$$
\end{lema}\qed

The result \cite[Th. 2.2]{DCY} can be extended to the comodule
algebra setting in a natural way. Let $U$ be a lifting of $H$
corresponding to the $l$-deformation $H_l$ under Theorem
\ref{main-dyc}.

\begin{teo} There is a natural bijection $Iso^{\infty}(G,H_\infty)\simeq
Lift(G,U)$.
\end{teo}

The following technical Lemma will be crucial to find liftings of
certain comodule algebras over pointed Hopf algebras.

\smallbreak

Let $F$ be a finite Abelian group, $H=\oplus_{i=0}^{m-1} H(i)$ be
a finite-dimensional coradically graded Hopf algebra,
$U^0\subseteq U^1\subseteq\dots \subseteq U^{m-1}=U$ a filtered
Hopf algebra $U$ with $U^0=\ku F=H(0)$, such that $\gr U=H$. Let
$g\in F$ and $x\in U^1-U^0$ an element such that
\begin{align}\label{commut-x}\Delta(x)=x\ot 1 + g\ot x,\quad
fx=\chi(f)\, xf,
\end{align}
for any $f\in F$ where $\chi$ is a character for $F$.

\smallbreak

Let $(G,\lambda_0)$ be a Loewy-graded $H$-comodule with grading
$G=\oplus_{i=0}^{m-1} G(i)$ such that $G(0)\simeq\ku
\widetilde{F}$ as $H$-comodules, where $\widetilde{F}\subseteq F$
is a subgroup, that is there exists a base $\{e_f:f\in
\widetilde{F}\}$ of $G$ such that the comodule structure are given
by $\lambda(e_f)=f\ot e_f,$ for any $f\in \widetilde{F}$.

\smallbreak

Let $(A, \lambda)$ be a left $U$-comodule algebra and let us
consider the filtration $A^i=\lambda^{-1}(U^i\otk A).$ Assume that
$\gr A=G$.

\begin{lema}\label{lambda1} Under the above assumptions,
if $\overline{y}\in G(1)$ is an element such that
$$ \lambda_0(\overline{y})= x\ot 1 + g\ot \overline{y},\quad
e_f \overline{y}= \chi(f)\; \overline{y} e_f,$$ for all $f\in
\widetilde{F}$, then there exists an element $y\in A^1$ such that
$\lambda(y)=x\ot 1 + g\ot y$, $e_f y= \chi(f)\; y e_f$ for all
$f\in \widetilde{F}$ and the class of $y$ in $A^1/A^0=G(1)$ equals
to $\overline{y}$.
\end{lema}
\pf Since we know that there are homogeneous maps $\lambda_s:G\to
H\otk G$ of degree $-s$, such that the coaction of $A$ is
$\lambda=\sum_{s=0}  \lambda_s$, then $\lambda(\overline{y})= x\ot
1 + g\ot \overline{y} + \lambda_1(\overline{y})$ where $
\lambda_1(\overline{y})\in H(0)\ot G(0)$. Let us write
$\omega=\lambda_1(\overline{y})$. Set $\omega=\sum_{h\in F, f\in
\widetilde{F}} \alpha_{h,f}\; h\ot e_f$, where $\alpha_{h,f}\in
\ku$. Since $(\id\ot\lambda)\lambda= (\Delta\ot\id)\lambda$ then
$$\sum_{h\in F, f\in \widetilde{F}} \alpha_{h,f}\; h\ot h\ot e_f=
\sum_{h\in F, f\in \widetilde{F}} \alpha_{h,f}\; h\ot f\ot e_f +
\sum_{h\in F, f\in \widetilde{F}} \alpha_{h,f}\; g\ot h\ot e_f,$$
whence we deduce that $\omega=\sum_{f\in \widetilde{F}, f\neq g}
\beta_f\; (g-f)\ot e_f$ for some $\beta_f\in \ku$. If we define
$a= \sum_{f\in \widetilde{F}, f\neq g}\beta_f\;  e_f$ then
$\omega= g\ot a-\lambda(a)$, thus the element $y=\overline{y}+a$
satisfies that $\lambda(y)=x\ot 1 + g\ot y$.

\smallbreak

The space $\mathcal{P}=\{y\in A: \lambda(y)=\mu\, x\ot 1 + g\ot y,
\mu\in \ku\}$ is stable under the action of $\widetilde{F}$ given
by
$$ f\cdot y= e_f y e_{f^{-1}}.$$
Thus there exists an element $y\in \mathcal{P}$ such that for any
$f\in \widetilde{F}$ $f\cdot y= \chi'(f)\, y$, for some
$\chi':\widetilde{F}\to \ku$. Since for all $f\in \widetilde{F}$
$fx=\chi(f)\, xf$ then $\chi'=\chi$ and this ends the proof of the
Lemma.

 \epf

\section{Deformation bicomplex for comodule
algebras}\label{section:bicomplex}

In \cite{y1} and \cite{Y2} the author introduced a bicomplex to
study algebraic deformations of $H$-comodule algebras $A$ over a
graded Hopf algebra $H$. The cohomology groups of the total
complex are related to formal deformations of $A$ over the same
Hopf algebra $H$. In this section we introduce a certain bicomplex
generalizing the one constructed by D. Yau in a way that
cohomology groups of the total complex are related to formal
deformations of $A$ over a any lifting $U$ of $H$.

\subsection{Bialgebra cohomology}

Let $H$ be a bialgebra. For each $s\in \Na$ and each $1\leq i\leq
s$ define the maps $m^s_i, \lambda^s_l, \lambda^s_r: H^{s+1}\to
H^s$ and $\Delta^s_i, \delta^s_l, \delta^s_r:H^s\to H^{s+1}$ by
\begin{align*} m^s_i(y^1\ot\dots \ot y^{s+1})&=y^1\ot\dots\ot y^i y^{i+1}\ot\dots \ot
y^{s+1},\\
\Delta^s_i(y^1\ot\dots y^s)&=y^1\ot\dots\ot y^i\_1\ot y^i_2 \ot \dots \ot y^s, \\
\lambda^s_l(x\ot y^1\ot\dots y^s)&= x\_1y^1\ot \dots\ot
x\_s y^s,\\
\lambda^s_r(y^1\ot\dots y^s\ot x)&=y^1x\_1\ot \dots  \ot y^s x\_s,\\
\delta^s_l(y^1\ot\dots y^s)&=y^1\_1\dots  y^s\_1\ot y^1\_2\ot
\dots \ot
y^s\_2 ,\\
\delta^s_r(y^1\ot\dots y^s)&=y^1\_1\ot\dots \ot y^s\_1\ot y^1\_2
\dots  y^s\_2,
\end{align*}
for all $x, y^1,\dots, y^s\in H$.

Let $C^{p,q}(H)=\Homk(H^{\ot p}, H^{\ot q})$ and let
$$\partial^h_{p,q}:C^{p,q} \to  C^{p+1,q}, \quad
\partial^v_{p,q}:C^{p,q} \to  C^{p,q+1},$$
be the maps defined by $\partial^h_{p,q}=\sum^{p+1}_{i=0}\,
(-1)^i\, \partial^h_{p,q}[i],$ and
$\partial^v_{p,q}=\sum^{q+1}_{i=0}(-1)^i\, \partial^v_{p,q}[i]$,
where for any $f\in C^{p,q}(H)$
$$\partial^v_{p,q}[i](f)=\begin{cases} (\id_H\ot f)\delta^p_l \;
&\text{ if }   i=0 \\
\Delta^q_i \circ f \;&\text{ if }  1\leq i\leq q\\
(f\ot\id_H)\delta^p_r \;&\text{ if } i=q+1
\end{cases} $$
and

$$\partial^h_{p,q}[i](f)=\begin{cases}\lambda^q_l (\id_H\ot f) \;
&\text{ if }   i=0 \\
 f\circ m^p_i \;&\text{ if }  1\leq i\leq p\\
\lambda^q_r(f\ot\id_H) \;&\text{ if } i=p+1
\end{cases} $$

We will use the truncated complex, that is $C^{p,q}(H)=0$ if
$pq=0$.

\subsection{Deformation bicomplex for comodule algebras}

Let $(A, \mu, \lambda)$ be a left $H$-comodule. For any $n\in \Na$
and $a^i\in A$ define the maps
$$\lambda^n_1:A^{\ot n}\to H\otk  A^{\ot n}, \quad
\lambda^n_2:A^{\ot n}\to H^{\ot n}\otk  A$$ by
$$\lambda^n_1(a^1\ot\dots\ot a^n)=a^1\_{-1}\dots  a^n\_{-1}\ot a^1\_0\ot\dots \ot
a^n\_0, $$
$$\lambda^n_2(a^1\ot\dots\ot a^n)=a^1\_{-1}\ot\dots  \ot a^n\_{-1}\ot a^1\_0\dots
a^n\_0, $$ for all $a^1\dots a^n\in A$. Define also an
$A$-bimodule structure on $H^{\ot n}\otk A$ via the maps
$$\beta^n_l:A\otk H^{\ot n}\otk A \to H^{\ot n}\otk A, \;\;\;\;
\beta^n_r:H^{\ot n}\otk A \otk A \to H^{\ot n}\otk A,$$ by
$$\beta^n_l(a\ot y^1\ot\dots\ot y^n\ot
b)=\lambda^n_l(a\_{-1}\ot y^1\ot\dots\ot y^n)\ot a\_0b,$$
$$\beta^n_r(y^1\ot\dots\ot y^n\ot a\ot b)=\lambda^n_r(y^1\ot\dots\ot y^n\ot
b\_{-1})\ot  ab\_0,$$ for all $y^1,\dots, y^n\in H$, $a,b\in A$.
For $n=0$ we define $\beta^0_l=\mu=\beta^0_r$. We also shall need
for any $n\in \Na$ and $1\leq i\leq n$ the maps
$$ \mu^n_i:A^{\ot n+1 }\to A^{\ot n+1 }, \quad
\mu^n_i(a^1\ot\dots\ot a^{n+1})=a^1\ot\dots\ot
a^ia^{i+1}\ot\dots\ot a^{n+1},$$ for all $a^1,\dots, a^{n+1}\in
A$.

Let us recall the bicomplex $C^{p,q}(A)$ defined by
$$C^{p,q}(A)=\begin{cases} 0  \quad &\text { if }   p=0\\
\Der(A) \quad &\text { if }  (p,q)=(1,0)\\
\Homk(A^{\ot p }, H^{\ot q }\otk A)\quad &\text { otherwise. }
\end{cases}$$
Define the horizontal and vertical differentials
$$\partial^h_{p,q}:C^{p,q}(A)\to C^{p+1,q}(A), \;\;\;
\partial^v_{p,q}:C^{p,q}(A)\to C^{p,q+1}(A),  $$
by
$$\partial^h_{p,q}=\sum^{p+1}_{i=0}(-1)^i\, \partial^h_{p,q}[i],\quad
 \partial^v_{p,q}=\sum^{q+1}_{i=0}(-1)^i\, \partial^v_{p,q}[i].$$
Where for any $f\in C^{p,q}(A)$
$$\partial^h_{p,q}[i](f)=\begin{cases} \beta^q_l (id_A\ot f)  &\text { if }
i=0\\
f\circ \mu^p_i &\text { if } 1\leq i\leq p\\

\beta^q_r (f\ot\id_A) &\text { if }  i=p+1,
\end{cases}$$
and

$$\partial^v_{p,q}[i](f)=\begin{cases} (\id_H\ot f) \lambda^p_1 &\text { if }
i=0\\
(\Delta^q_i\ot\id_A) f &\text { if } 1\leq i\leq q\\
(\id_{H^{\ot q}}\ot\lambda) f &\text { if }  i=q+1.
\end{cases}$$

We now define a new complex. For each $p,q$ set
$\ca^{p,q}(A,H)=C^{p,q}(A)\oplus C^{p,q+1}(H)$. Define the
horizontal and vertical differentials by
$$d^h_{p,q}:\ca^{p,q}(A,H)\to \ca^{p+1,q}(A,H), \quad
d^v_{p,q}:\ca^{p,q}(A,H)\to \ca^{p,q+1}(A,H)$$ by the formulas
$$d^h_{p,q}(\phi, f)=(\partial^h_{p,q}\phi\; ,\, \partial^h_{p,q+1}
f),$$ and $$d^v_{p,q}(\phi, f)=(\partial^v_{p,q}\phi + (-1)^{q}\;
(f\ot\id_A)\lambda^p_2\; ,\, \partial^v_{p,q+1}f).$$ Here we abuse
of the notation by using the same symbol for denoting the
differentials coming from the bialgebra cohomology of $H$ and
those coming from the cohomology of the comodule algebra $A$.

\begin{lema} $\ca^{*,*}(A,H)=(\ca^{p,q}(A,H),d^h_{p,q}, (-1)^{p} d^v_{p,q}) $ is a bicomplex.
\end{lema}
\pf We have to verify that the following equations hold: $d^h\circ
d^h=0$, $d^v\circ d^v=0$ and $d^h\circ d^v=d^v\circ d^h$. The
first one is obvious. Let us prove the other two. \medbreak

It can be easily checked, case by case, that for any $0\leq i\leq
q+2$ and any $0\leq j\leq p+1$ the equalities $$\partial^v_{p,q+1
}[i]((f\ot\id_A)\lambda^p_2)=
(\partial^v_{p,q}[i](f)\ot\id_A)\lambda^p_2, $$
$$\partial^h_{p,q+1 }[j]\big((f\ot\id_A) \lambda^p_2\big)=
\big(\partial^h_{p,q+1 }[j](f)\ot\id \big) \lambda^{p+1}_2$$ hold.
This implies that \begin{equation}\label{d2=01} \partial^v_{p,q+1
}((f\ot\id_A)\lambda^p_2)=
(\partial^v_{p,q}(f)\ot\id_A)\lambda^p_2\end{equation} and
\begin{equation}\label{d2=02} \partial^h_{p,q+1 }\big((f\ot\id_A) \lambda^p_2\big)=
\big(\partial^h_{p,q+1 }(f)\ot\id \big)
\lambda^{p+1}_2.\end{equation}

\medbreak

Let $p,q\in \Na$ and $(\phi,f)\in \ca^{p,q}(A,H)$ then
\begin{align*} d^v_{p,q+1} d^v_{p,q}(\phi,f)&=
d^v_{p,q+1} \big( \partial^v_{p,q}\phi + (-1)^{q}\;
(f\ot\id_A)\lambda^p_2\; ,\, \partial^v_{p,q}f\big)\\
&=\big(\partial^v_{p,q+1 }\partial^v_{p,q}\phi +
(-1)^{q}\partial^v_{p,q+1 }((f\ot\id_A)\lambda^p_2)+\\&
(-1)^{q+1}(\partial^v_{p,q}f\ot\id_A)\lambda^p_2,
\partial^v_{p,q+1}\partial^v_{p,q}f\big).
\end{align*}
Follows from \eqref{d2=01} that this last expression is zero. By a
straightforward computation one proves that $d^v_{p+1,q}
d^h_{p,q}=d^h_{p,q+1}d^v_{p,q}$ using \eqref{d2=02}.

\epf

The total complex $(\Tot^*(\ca^{*,*}(A,H)), d)$ of the bicomplex
$\ca^{*,*}(A,H)$ is defined by
$$\Tot^n(\ca^{*,*}(A,H))=\bigoplus_{p+q=n} C^{p,q}(A) \oplus
\bigoplus_{p+q=n}  C^{p,q+1}(H),$$ and the differentials
$d^n:\Tot^n(\ca^{*,*}(A,H))\to \Tot^{n+1}(\ca^{*,*}(A,H))$ are
defined by
$$d^n|_{\ca^{n-i,i}(A,H)}= d^h_{n-i,i} + (-1)^p\, d^v_{n-i,i},$$
for each $0\leq i\leq n$. Borrowing the notation of D. Yau
\cite{y1}, we will denote the cohomology groups of this complex by
$H^*_{ca}(A,H)$.

\medbreak

Let $(\phi, f)\in \ca^{2,0}(A,H)=C^{2,0}(A)\oplus C^{2,1}(H)$ ,
$(\psi, g)\in  \ca^{1,1}(A,H)=C^{1,1}(A)\oplus C^{1,2}(H)$. The
element $(\phi, f)+ (\psi, g)\in Z^2_{ca}(A,H)$ if and only if
\begin{align} \partial^h_{2,0}\phi=0,\quad \partial^v_{1,1}
\psi-(g\ot\id_A)\lambda=0
\end{align}
\begin{align}\partial^h_{2,1} f=0=\partial^v_{1,2} g, \quad
\partial^v_{2,1} f+ \partial^h_{1,2} g=0,
\end{align}
\begin{align}\partial^v_{2,0}\phi + (f\ot\id_A)\lambda^2_2 +
\partial^h_{1,1}\psi=0.
\end{align}

The next result follows from a simple computation.
\begin{prop} Let $H$ be a graded Hopf algebra, $(G,\mu,\lambda)$
a graded left $H$-comodule algebra. Let $\big(H[t]/ (t^{l+1}),
m^l, \Delta^l\big)$ be an $l$-deformation of $H$ and $\big(G[t]/
(t^{l+1}), \mu^l, \lambda^l\big)$ an $l$-deformation of $G$
compatible with the deformation of $H$. Then $(\mu^l_1+m^l_1) +
(\lambda^l_1, \Delta^l_1)\in Z^2_{ca}(G,H)$.\qed
\end{prop}

\section{Comodule algebras over pointed Hopf
algebras}\label{section:comod-pointed}

Let $H$ be a Hopf algebra and $R$ be a braided Hopf algebra in $
{}^{H}_{ H} \mathcal{YD}$.

\begin{defi} A \emph{braided left comodule algebra} over $R$ is an
algebra $B\in  {}^{H}_{ H} \mathcal{YD}$ together with a linear
map $\delta_B: B\to R\otk B$, $\delta_B(x)=x^{[-1]}\ot x^{[0]}$,
such that $\delta_B$ is a morphism in the category ${}^{H}_{ H}
\mathcal{YD}$, $(B, \delta_B)$ is a left $R$-comodule, and such
that for all $x, y\in B$
\begin{equation} \delta_B(xy)= x^{[-1]} (x^{[0]}\_{-1}\cdot\,
y^{[-1]})\ot x^{[0]}\_0 y^{[0]}.
\end{equation}
\end{defi}

\subsection{Comodule algebras over bosonizations}

We shall give a simple recipe to produce comodule algebras over a
Hopf algebra constructed from a bosonization.

\medbreak

Let $H_0$ be a Hopf algebra and $R\in  {}^{H_0}_{ H_0}
\mathcal{YD}$ be a braided Hopf algebra. Let us recall that the
bosonization, or Radford biproduct, $H=R\# H_0$, is the Hopf
algebra over the vector space $R\otk H_0$ with product and
coproduct given by
$$(r\# g)(s \# f)= r g\_1\cdot s \# g\_2 f, $$
$$\Delta( r\# g)= r\_1\#(r\_2)\_{-1} g\_1 \ot (r\_2)\_0\#g\_2,$$
for all $r, s\in R$, $f, g\in H$.

\medbreak

Let $B$ be a braided comodule algebra over $R$ and $F$ be a left
$H_0$-comodule algebra. We shall give to the tensor product $B\otk
F$ a natural structure of $R\# H_0$-comodule algebra. Let us
define the product and left $R\# H_0$-comodule structure by:
$$\delta(x\ot f)=  x^{[-1]}\# \, x^{[0]}\_{-1} f\_ {-1}\ot\, x^{[0]}\_0\ot f\_0,$$
$$ (x\ot f)(y\ot g)= x (f\_{-1}\cdot y)\ot f\_0 g,$$
for all $x,y \in B$, $f,g\in F$. We shall denote by $B\# F$ the
space $B\otk F$ together this new product and the coaction
$\delta$.

\begin{prop}\label{smash-comodule-algebra} The following
assertions hold:
\begin{enumerate}
    \item $B\# F$ is a left $R\# H_0$-comodule algebra.
    \item If $F^{\co H_0}=\ku $ and $B^{\co R}=\ku$ then we have
 that $(B\# F)^{\co R\# H_0}=\ku $.

\end{enumerate}
\end{prop}
\pf (1) follows by a straightforward computation.

Let $\sum x_i\ot f_i \in (B\# F)^{\co R\# H_0}$. Then
$$\delta(\sum x_i\ot
f_i)=\sum x_i^{[-1]}\# \, x_i^{[0]}\_{-1} f_i\_ {-1}\ot\,
x_i^{[0]}\_0\ot f_i\_0=\sum 1\ot 1 \ot x_i\ot f_i.
$$
Applying the counit on the second tensorand we obtain that
$$\sum x_i^{[-1]}\ot x_i^{[0]}\ot f_i= \sum 1 \ot x_i\ot f_i,$$
whence $ x_i\in B^{\co R}=\ku$, which implies that $f_i\in F^{\co
H_0}=\ku $. This proves (2).\epf

\subsection{Loewy-graded  comodule algebras over pointed Hopf
algebras}

For the rest of this section $G$ will denote a finite group, $V\in
{}^{\ku G}_{\ku G} \mathcal{YD} $ will denote a Yetter-Drinfeld
module over $G$ and $\nic(V)$ the associated Nichols algebra, see
\cite{AS2} and references therein. The Nichols algebra is a
$\Na_0$-graded braided Hopf algebra with grading
$\nic(V)=\oplus_{i=0}^m \nic^i(V)$. Let us denote by $\delta: V\to
\ku [G]\otk V$ the coaction.

\medbreak

We denote the bosonization  $H=\nic(V)\#\ku [G]$. The Hopf algebra
$H$ is a coradically graded Hopf algebra, with grading given by
$H(n)=\nic^n(V)\otk \ku [G]$ for any $n=0\dots m$. We will denote
by $\theta:H\to \nic(V)$ and $p:H\to \ku [G]$ the canonical
projections. For any $h\in H$ we have that $h=\theta(h\_1)
p(h\_2)$. The space $\nic(V)$ are the coinvariants of $H$ under
the coaction $(\id\ot p)\Delta$, that is $\nic(V)=\{h\in h:(\id\ot
p)\Delta(h) =h\ot 1\}$.

\medbreak

Let $(A,\lambda)$ be a Loewy-graded left $H$-comodule algebra with
gradation given by $A=\bigoplus^m_{i=0}\,A(i)$. Let $\pi:A\to
A(0)$ be the canonical projection.

Let us define
$$\nic_A=\{a\in A: (\id\ot\pi)\lambda(a)\in H\ot 1\}.$$
This space should be thought as a kind of \emph{diagram}
\cite{AS2}, for the comodule algebra $A$. Indeed, if
$(A,\lambda)=(H,\Delta)$ then $\nic_A=\nic(V)$ is the diagram
corresponding to $H$.

\medbreak

For each $n\in \Na$ define $\nic_A(n)=\nic_A\bigcap A(n)$.

\begin{prop}\label{diagram-comod} The following statements holds.
\begin{enumerate}\item $\nic_A(0)=\ku 1$,
$\nic_A= \bigoplus_{n=0} \nic_A(n)$.

    \item $\nic_A\subseteq A$ is a left $H$-subcomodule subalgebra,
    in particular $\lambda(\nic_A(n))\subseteq \oplus_{i=0}^n\;
    H(n-i)\otk \nic_A(i)$.
    \item There exists an injective $H$-comodule algebra
    homomorphism $\iota:\nic_A\to \nic(V)$ such that $\iota(\nic_A(n))\subseteq
    \nic^n(V)$ for any $n=0\dots m$.

\item If $\nic_A(1)=0$ then $\nic_A=\ku 1$.

 \item The space
$\iota(\nic_A(1))\subseteq V$ is a $\ku G$-subcomodule.

    \item The multiplication map $\mu:\nic_A\otk A_0\to A$ is injective.

\end{enumerate}
\end{prop}
\pf (1). The first assertion is clear. Let $a\in \nic_A$, then
$a=\sum_{i=0}\, a_i$, where $a_i\in A(i)$. Since $(\id\ot
\pi)\lambda(a)\in H\otk 1$ and the spaces $H(i)$ are disjoint we
conclude that $(\id\ot \pi)\lambda(a_i)\in H(i)\otk 1$ for each
$i$, hence $a_i\in \nic_A(i).$

\medbreak

(2). It is clear that $\nic_A$ is a subalgebra since $\lambda$ and
$\pi$ are algebra maps. Let $a\in \nic_A$, then $a\_{-1}\ot
\pi(a\_0)= h\ot 1$, for some $h\in H$. Then
\begin{equation}\label{iota-comod}\Delta(h)\ot 1=a\_{-2}\ot a\_{-1}\ot\pi(a\_0)=a\_{-1} \ot a\_0\_{-1}\ot
\pi(a\_0\_0).\end{equation} Therefore $ a\_{-1}\ot a\_0 \in H\ot
\nic_A$.

\medbreak

(3). For $a\in\nic_A$ define the map $\iota:\nic_A\to H$ by
$\iota(a)=h\in H$, where $h$ is the unique element in $H$ such
that $a\_{-1}\ot \pi(a\_0)= h\ot 1$. Using \eqref{iota-comod} it
is easy to see that $\iota$ is an $H$-comodule map.

Let us prove that $\iota$ is injective. Let $a\in\nic_A$ such that
$\iota(a)=0$, thus $a\_{-1}\ot \pi(a\_0)=0$.  Let us write
$a=\sum_n\; a_n$ where $a_n\in A(n)$. In this case
$\lambda(a)=\sum_n\lambda(a_n)=\sum_n\sum^n_{i=0}\, b^n_i$, where
$b^n_i\in H(i)\otk A(n-i)$, for $i=0\dots n$. Since
$a\_{-1}\ot\pi(a\_0)=\sum_n b^n_n=0$ this means that for each $n$,
$b^n_n=0$, thus  $\lambda(a_n)\in H_{n-1}\otk A$, whence $a_n\in
A_{n-1}$, which is impossible unless $a_n=0$ since for all $n$:
$A(n)\bigcap A_{n-1}=0$. Thus $\iota$ is injective. \medbreak

Observe that $\lambda\pi=(p\ot \pi) \lambda$, thus if $a\in
\nic_A$ then there exists an element $h\in H$ such that
$a\_{-1}\ot \pi(a\_0)= h\ot 1$. Hence
\begin{align*}(\id\ot p) \Delta)(h)\ot 1&= a\_{-2}\ot p(a\_{-1})\ot\pi(a\_0)=
a\_{-1}\ot (p\ot \pi)\lambda(a\_0)\\
&=a\_{-1}\ot\lambda\pi(a\_0)=h\ot 1\ot 1.
\end{align*}
Thus $(\id\ot p) \Delta(h)=h\ot 1$, that is $h\in \nic(V)$. Since
$H(n)= \nic^n(V)\# \ku G$ and $\iota$ is an $H$-comodule map then
$\iota(\nic_A(n))\subseteq
    \nic^n(V)$. This ends the proof of (3).

\medbreak

(4). Let us assume that $\nic_A(1)=0$. We can (and will) assume
that $\nic_A\subseteq \nic(V).$ Let $a\in\nic_A(2)$ then
$\Delta(a)\in H_0\otk \nic_A(2) \oplus H(2)\otk \nic_A(0)$, thus
$a\in H_1$, hence $a=0$. Using the same argument we can prove
inductively that $\nic_A(n)=0$ for all $n\geq 2$.

\medbreak

(5). Let us consider  $W= \nic_A(1)$ as a $\ku G$-comodule via
$(p\ot\id_A) \lambda$. Then $W=\oplus_{g\in G} W_g$, where
$W_g=\{w\in \nic_A(1): p(a\_{-1})\ot a_0= g\ot a\}$. If $a\in W_g$
then $\lambda(a)= g\ot a + v_a\ot 1$, for some $v_a\in V$. Note
that $\iota(a)=v_a$. In this case, $\delta(v_a)=g\ot v_a$, thus
$\iota(W)$ is a subcomodule of $V$.

\medbreak

Now, let us prove (6). Consider the map $\phi:A\to \nic(V)\otk
A_0$ given by $\phi(a)= \theta(a\_{-1})\ot \pi(a\_0)$ and let
$\mu:\nic_A\otk A_0\to A$ be the multiplication map. Then
$\phi\circ \mu= \iota\ot \id_{A_0}$, indeed if $x\in \nic_A, a\in
A_0$ then
\begin{align*}\phi(xa)&=\theta(x\_{-1}a\_{-1})\ot \pi(x\_{-1}a\_0) =
\theta(x\_{-1}a\_{-1})\ot \pi(x\_{-1}) a\_0\\
&=\theta(\iota(x)a\_{-1})\ot a\_0=\iota(x)\ot a.
\end{align*}
Since $\iota$ is injective this implies that $\mu$ is injective.

\epf

\begin{rmk} The space $\iota(\nic_A)$ is a left homogeneous coideal
subalgebra of $\nic(V)$. The algebra $\nic_A$ need not be
generated by $\nic_A(1)$. However, in some special cases, for
example when $\nic(V)$ is a quantum linear space, $\nic_A$ is
always generated in degree one.
\end{rmk}

Let us further assume that $(A,\lambda)$ is right $H$-simple.
Follows by Proposition \ref{equivalence-grad} that $A(0)$ is right
$\ku [G]$-simple and therefore there exists a subgroup $F\subseteq
G$ and a 2-cocycle $\psi\in Z^2(F,\ku^{\times})$ such that
$A(0)\simeq \ku^{\psi}[F]$ as $\ku [G]$-comodule algebras. That
is, there is a basis $\{e_f\}_{f\in F}$ of $A(0)$ such that
$$ e_f e_{f'}=\psi(f,f')\, e_{ff'},\quad \lambda(e_f)=f\ot e_f,$$
for all $f, f'\in F$.

\smallbreak

In the case when $\psi$ is trivial there is an action of $F$ on
$\nic_A$ given as follows. If $x\in \nic_A$ and $f\in F$ then
$$ f\cdot x:= e_f x e_{f^-1} \in \nic_A.$$
With this action and the trivial coaction $\delta:\nic_A\to
\ku[F]\otk \nic_A$, given by $\delta(x)= 1\ot x$, $\nic_A$ is a
Yetter-Drinfeld module over $\ku[F]$. Thus, we can consider the
smash product $\nic_A\# A(0)$ as in Proposition
\ref{smash-comodule-algebra}. \bigbreak

The following proposition tells us that the comodule algebra $A$
can be recovered from $\nic_A$ and $A_0$.
\begin{prop}\label{diagram-comod2} Assume that the group $G$ is a cyclic
group, then $A\simeq  \nic_A\# A(0)$ as $H$-comodule algebras.
\end{prop}

\pf Since $G$ is cyclic then $\psi$ must be trivial and $A_0$ is
isomorphic to the group algebra $\ku F$, where $F\subseteq G$ is a
subgroup. Thus, we can consider the product $ \nic_A\# A(0)$.

By Proposition \ref{diagram-comod} (6) the multiplication map
$\mu:\nic_A\# A_0\to A$ is injective and clearly is an
$H$-comodule algebra map. Let us prove that $\mu$ is surjective.

\medbreak

Let $Q$ denote the quotient coalgebra $H/(\ku F)^+ H$ and
$\gamma:H\to Q$ the canonical projection. Since $A\in {}^{H}_{\ku
F}\Mo$, then $A\simeq H\square_Q \overline{A}$, where
$\overline{A}=A/(\ku F)^+ A$ and the left $Q$-comodule structure
on $\overline{A}$ is given by $(\gamma\ot\id_A)\lambda$. The left
$H$-comodule structure on $H\square_Q \overline{A}$ is given by
the coproduct of $H$.

\medbreak

There is an isomorphism $\alpha:\ku F \ot
Q\xrightarrow{\;\;\simeq\;\;} H $
 of left $\ku F$-modules, right
$Q$-comodules. Thus there is a linear isomorphism $A\simeq \ku F
\ot\overline{A}$. \smallbreak

It is easy to see that $(H\square_Q \overline{A})_0= \ku F\ot
\overline{1}$. Let us denote by $\pi:H\square_Q \overline{A}\to
(H\square_Q \overline{A})_0$ the projection. The image of the
injective map $\beta:\overline{A}\to H\square_Q \overline{A}$
given by $\beta(\overline{a})= \alpha(1\ot\overline{a}\_{-1}) \ot
\overline{a}\_0$ composed with the isomorphism $A\simeq H\square_Q
\overline{A}$ is inside $\nic_A$. Indeed, for any $\overline{a}\in
\overline{A}$ the element
$$ \alpha(1\ot\overline{a}\_{-1})\_1 \ot\; \pi( \alpha(1\ot\overline{a}\_{-1})\_2
\ot \overline{a}\_0) \in H\ot 1.$$ This implies that $\dim
\overline{A}\leq \dim \nic_A$. Whence, $\dim A \leq \nic_A \mid
F\mid$, thus $\mu$ is surjective and therefore a bijection.\epf

\begin{rmk}\label{stable-F} Assume we are under the hypothesis of
Proposition \ref{diagram-comod2}. It is immediate to verify that
under the identification of $A(0)\simeq \ku F \hookrightarrow \ku
G$, the space $\iota(\nic_A)\subseteq \nic(V)$ is stable under the
action of $F$.
\end{rmk}

As a remarkable consequence of the above proposition we have the
following result. Let $U$ be a Hopf algebra such that $\grc
U=\nic(V)\# \ku G$ for a cyclic group $G$.

\begin{lema} If $A$ is a $U$-simple left $U$-comodule algebra with
trivial coinvariants, then $\dim A$ divides $\dim U$.\qed
\end{lema}

This result is related with a conjecture made by Y. Zhu for
semisimple Hopf algebras and it is related to Kaplansky's sixth
conjecture. See \cite[Conjecture 4.1]{Z}.

\section{Some classification
results}\label{section:classification}

Let $H$ be a pointed Hopf algebra with coradical the group algebra
of a cyclic group $C_n$, such that $\grc H=\nic(V)\# \ku[C_n]$,
where $V\in {}_{\ku C_n}^{\ku C_n}\YD$ is a Yetter-Drinfeld module
and $\nic(V)$ the associated Nichols algebra.

\medbreak

 We propose the following strategy to classify module
categories over the category $\Rep(H)$. First we determine  all
Loewy-graded right $\nic(V)\# \ku[G]$-simple left $\nic(V)\#
\ku[G]$-comodule algebras $A$, with trivial coinvariants. By
Proposition \ref{diagram-comod2} this comodule algebras are
determined by the subalgebra $\nic_A$ and $A_0$, where
$\nic_A=\oplus_i \nic_A(i)$ is a homogeneous left coideal
subalgebra of $\nic(V)$ stable under the action of $A_0$. Since
$A_0$ is a $\ku C_n$-simple left $\ku C_n$-comodule algebra, then
$A_0=\ku C_d$ for some divisor $d$ of $n$. In our examples
$\nic_A$ is generated as an algebra by $\nic_A(1)$, so really
$\nic_A$ depends on a subspace $W\subseteq V$ stable under the
action of $C_d$. To end we shall find all possible liftings of
this comodule algebras.

\medbreak

In this section we present the classification of exact module
categories over some Hopf algebras, where the dimension of $V$ is
1 or 2. More precisely we will show the classification for the
Taft Hopf algebras $T_q$, over the Radford Hopf algebras
$\textbf{r}_q$, over the book Hopf algebras and over the small
quantum groups $u_q(\mathfrak{sl}_2)$.

\smallbreak
 Let us recall the definition of these Hopf algebras and some
other Hopf algebras that will be used later. Some notations are
taken from \cite{AS1}. Let $n\in \Na$ and $q$ be a $n$-th
primitive root of unity.

\medbreak
\begin{itemize}
    \item The \emph{Taft Hopf algebra}  $T_q= \ku\langle g, x\vert \; gx = q\,xg,
\;g^{n} = 1, \; x^n =0 \rangle,$ with coproduct  determined by
$\Delta(g) = g\otimes g$,
 $\Delta(x) = x \otimes 1+ g\otimes x $.\smallbreak
    \item The algebra $\widehat{T(q)} = \ku\langle g, x\vert  \; gx =
q\,xg, \;g^{n^2} = 1,  \;x^n = 0\rangle$. The coproduct is
determined by $\Delta(x) = x\otimes 1 + g \otimes x$, $\Delta(g) =
g\otimes g$.\smallbreak
    \item The \emph{Radford Hopf algebra} is $\bold r(q) = \ku\langle g,
x\vert \; gx = q\,xg, \;g^{n^2} = 1, \; x^n = 1 - g^n\rangle.$
 The coproduct is determined by $\Delta(g) = g\otimes g$,
 $\Delta(x) = g\otimes x + x \otimes 1$.\smallbreak

\item The \emph{book Hopf algebras}
 $\Hc(1,q)= \ku\langle g, x, y\vert \; g^n=1, \;
gx = q\, xg,  \; gy = q^{-1}\,yg,\; xy=q\, yx, \;
x^n=0=y^n\rangle.$ With coproduct determined by $\Delta(g) =
g\otimes g$, $\Delta(x) = x\otimes 1 + g^{-1} \otimes x$,
$\Delta(y) = y\otimes 1 + g^{-1} \otimes y$.

\item The \emph{Frobenius-Lusztig kernel} $u_q(\mathfrak{sl}_2) =
\ku\langle g, x, y\vert  \; gx = q^2\, xg,  \; gy = q^{-2}\,yg,
 \;g^{n} = 1,\; x^n = 0= y^n,  \; xy-yx =
 g-g^{-1}\rangle.$
The coproduct is determined by $\Delta(g) = g\otimes g$,
$\Delta(x) = x\otimes g + 1 \otimes x$, $\Delta(y) = y\otimes 1 +
g^{-1} \otimes y$.

\end{itemize}

\medbreak

\subsection{Module categories over $\Rep(T_q)$ } In this
section we recover, using different techniques, the results
obtained in \cite[Thm 4.10]{eo}.

\medbreak

Let $d$ be a divisor of $n$. Set $n=dm$.  For any $\xi\in\ku$
define the algebra $\Ac(d,\xi)$ generated by elements $h, w$
subject to relations
$$ h^d=1,\;\;\; hw=q^m\; wh,\;\;\; y^n=\xi\, 1.$$
The algebra $\Ac(d,\xi)$ is a left $T_q$-comodule algebra with
coaction determined by
$$\lambda(h)=g^m\ot h, \; \lambda(w)= x\ot 1 + g\ot w.$$

Let us denote by $K(d)$ the left coideal subalgebra of $T_q$
generated by $g^m$ and $x$.

\begin{lema}\label{tec21}For any $\xi\in\ku$ there is an isomorphism
$\Ac(d,\xi)\simeq K(d)_{\sigma_\xi}$ for some  cocentral 2-cocycle
$\sigma_\xi$. In particular the algebras $\Ac(d,\xi)$ are right
$T_q$-simple left $T_q$-comodule algebras such that
$\Ac(d,\xi)^{\co T_q}=\ku$.
\end{lema}

\pf The algebra $\Ac(n,\xi)$ is a left $T_q$-Galois extension and
$\Ac(d,\xi)\subset \Ac(n,\xi)$ is a subalgebra and a left
$T_q$-submodule. Thus there exists a 2-cocycle such that
$\Ac(n,\xi)\simeq {}_{\sigma_\xi} T_q$. This cocycle is cocentral
since $\Ac(n,\xi)$ is also a right $T_q$-comodule algebra with
structure map $\rho:\Ac(n,\xi)\to \Ac(n,\xi)\otk T_q$ given by
$$ \rho(h)=h\ot g, \quad \rho(w)=w\ot 1 + h\ot x.$$
The other statements are straightforward.\epf

\begin{lema}

\begin{itemize}

   \item[(i)]   Categories ${}_{\Ac(d,\xi)}\Mo$
    are semisimple for any divisor $d$ of $n$ and any $\xi\in
    \ku^{\times}$.

    \item[(ii)] For any $\xi\in
    \ku^{\times}$ the rank of ${}_{\Ac(d,\xi)}\Mo$ is
    $\frac{n}{d}.$
    \item[(iii)] The algebras $\Ac(d,\xi)$, $\Ac(d',\xi')$
     are Morita equivariant equivalent if and only if
     $d=d'$, $\xi=\xi'$.
\end{itemize}
\end{lema}
\pf (i) and (ii)  are Straightforward. Let us assume that
$\Ac(d,\xi)\sim_M \Ac(d,\xi')$. From (ii) we get that $d=d'$. Let
us denote by $Q$ the quotient coalgebra $Q=T_q/K(d)^+T_q\simeq \ku
C_m$. Using Lemma \ref{morita-cocentral} and Lemma \ref{tec21}
there exists a group-like element $f\in G(T_q)$ such that
$\Ac(d,\xi') \simeq (f K(d)f^{-1})_{\sigma_\xi}$. Thus
$\Ac(d,\xi')\simeq \Ac(d,\xi)$, whence $\xi=\xi'.$\epf

\begin{prop}\label{modcat-taft} If $\Mo$ is an exact indecomposable module category
over $\Rep(T_q)$ then $\Mo\simeq {}_{\ku C_d}\Mo$ or $\Mo\simeq
{}_{\Ac(d,\xi)}\Mo$ for some divisor $d$ of $n$ and $\xi\in \ku$.

\end{prop}

\pf Let $G$ be a Loewy-graded left $T_q$-comodule algebra. Let $d$
be a divisor of $n$ such that $G_0=\ku C_d$. Let us assume that
$G\neq G_0$. Using Proposition \ref{diagram-comod2} $G\simeq
\nic_G\# \ku C_d$. Since $G\neq G_0$ then $\nic_G(1)\neq 0$. The
only possibility is that $G=\Ac(d,0)$.

\medbreak

Let $(A, \lambda)$ be a  $T_q$-simple left $T_q$-comodule algebra
such that $A$ is a lifting of $\Ac(d,0)$ along $T_q$ for some
divisor $d$ of $n$. So in this case $A_0=\ku C_d$ is the group
algebra of the cyclic group of $d$ elements generated by $h$ where
$\lambda(h^i)=g^{mi}\ot h^i$, for any $i=0\dots d-1$, where
$dm=n$.

By Lemma \ref{lambda1} there is an element $y\in \Ac_1- \Ac_0$
such that $\lambda(y)=x\ot 1 + g\ot y$ and and $hy=q^m\, yh$.
Since $A$ is a lifting for $\Ac(d,0)$, the set $\{h^iy^j: i=0\dots
d-1, j=0\dots n-1 \}$ is a basis for $A$. Here
$A_0=\ku\{h^i:i=0\dots d-1\}$.

\medbreak

Since $\lambda(y^n)=x^n\ot 1 + g^n\ot y^n= 1\ot y^n$ then $y^n\in
A_0$ and there exists an $\xi\in\ku$ such that $y^n=\xi 1$ and we
have a projection $\Ac(d,\xi)\twoheadrightarrow A$ which is an
isomorphism because both algebras have the same dimension. \epf

\subsection{Module categories over $\Rep(\textbf{r}_q)$ }

Let $C_{n^2}$ be the cyclic group of $n^2$ elements generated by
$g$. Let $V\in {}^{\ku C_{n^2}}_{\ku C_{n^2}} \mathcal{YD}$ be the
one-dimensional module generated by $w$ with action and coaction
determined by
$$g\cdot w= q\; w, \quad \delta(w)= g\ot w. $$

The following result is straightforward.

\begin{lema}\label{alg-grad} $\grc \bold r(q)=\widehat{T(q)}$ and
$\widehat{T(q)}\simeq \nic(V)\# \ku C_{n^2}$.\qed
\end{lema}

We shall classify exact module categories over the Radford Hopf
algebras in a similar way as for the Taft Hopf algebras. No new
difficulty arise since the dimension of the the vector space $V$
is 1. \medbreak

Let us define  families of left $\bold r(q)$-module algebras,
right $\bold r(q)$-simple. Let $d$ be a divisor of $n^2$ and set
$n^2=dm$. Let also $\xi\in \ku^{\times}$.

\begin{itemize}
    \item[(a)] The group algebra $\ku C_d=  \ku\langle h:
    h^d=1\rangle$, with coaction determined by $\lambda(h)= g^m\ot
    h$.
\item[(b)] Algebras $\Ac(d)=\ku\langle h, w:  h^d=1,\, w^n=1,\,
hw= q^m\; wh\rangle$ and coaction determined by $\lambda(h)=g^m\ot
    h$, $\lambda(w)= x \otimes 1 + g\otimes w$.

    \item[(c)] If there is an integer $a$ such that $na=d$,
     then $\Bc(a,\xi)=\ku\langle h, w:  h^d=1,\, w^n=1+ \xi\, h^a,\, hw=
    q^m\; wh\rangle$ and coaction determined by $\lambda(h)=g^m\ot
    h$, $\lambda(w)= x \otimes 1 + g\otimes w$.
\end{itemize}

For any $d$  divisor of $n^2$ and any $\xi\in \ku$ the algebras
$\Ac(d)$ and $\Bc(a,\xi)$ are left $\bold r(q)$-module algebras,
right $\bold r(q)$-simple with trivial coinvariants. The algebras
listed above are non-isomorphic as comodule algebras.

\begin{lema}\label{non-equivalent-radford} The algebras $\Ac(d)$ and $\Bc(a,\xi)$ are twistings
of left coideal subalgebras of $\textbf{r}(q)$ by a compatible
Hopf 2-cocycle, pair-wise non Morita equivariant equivalent.
\end{lema}
\pf Algebras $\Ac(n)$ and $\Bc(n,\xi)$ are $\textbf{r}(q)$-Galois
extensions of the field and $\Ac(d)\subseteq \Ac(n)$,
$\Bc(a,\xi)\subseteq \Bc(n,\xi)$ are $\textbf{r}(q)$-subcomodules
subalgebras. Let us prove now that these algebras are non Morita
equivariant equivalent. For this we will need the following
result:

\begin{claim}\label{P-in-grad} Let $H$ be a pointed Hopf algebra.
Let  $J\subseteq H$ be a left coideal subalgebra and $\sigma:H\otk
H\to \ku$ a compatible 2-cocycle, such that the graded algebra
$\grl \big({}_{\sigma}J\big)$ is isomorphic to a left coideal
subalgebra $K$ of $\grc H$ and the coalgebra quotient $Q=\grc
H/(\grc H) K^+$ is pointed cosemisimple. Then if $P\in
{}^H\Mo_{{}_{\sigma}J}$ is an indecomposable object there exists a
group-like element $g\in H_0$ such that $P\simeq g\cdot
{}_{\sigma}J$.
\end{claim}
\pf[Idea of the proof] Consider the filtration
$P_i=\delta^{-1}(H_i\otk P)$. The associated graded vector space
$\gr P$ is an object in the category ${}^{\grc H}\Mo_{K}\simeq
{}^Q\Mo$. Since $P$ is indecomposable, the so is $\gr P$.
Therefore there exists a group-like element $g\in H$ such that
$\gr P=gK$. \epf

Let $A$ be a left $\bold r(q)$-comodule algebra. Let us assume,
for instance, that $A\sim_M \Bc(a,\xi)$ for some $a$ such that
$na=d$. Then there exists an object $P\in
{}^{r(q)}\Mo_{\Bc(a,\xi)}$ such that $A\simeq
\End_{\Bc(a,\xi)}(P)$. Since in our case we are under the
hypothesis of Claim \ref{P-in-grad}, $P=g\cdot \Bc(a,\xi)$, and
arguing as in Lemma \ref{morita-cocentral} we obtain that $A\simeq
\Bc(a,\xi)$, thus if $A$ is an algebra in the list it must be
equal to $\Bc(a,\xi)$. The same argument can be used to prove that
if $A\sim_M \Ac(d)$ then $A\simeq \Ac(d)$.\epf

\begin{prop} If $\Mo$ is an exact indecomposable module category
over $\Rep(\textbf{r}(q))$ then $\Mo\simeq {}_{\ku C_m}\Mo$,
$\Mo\simeq  {}_{\Ac(d)}\Mo$ or $\Mo\simeq {}_{\Bc(a,\xi)}\Mo$.
\end{prop}

\pf Let $(G, \lambda_0)$ be a Loewy-graded
$\widehat{T(q)}$-comodule algebra. We assume that $G\neq G_0$ and
$G_0=\ku\langle h: h^d=1\rangle$ for some $d$ divisor of $n^2$.
Since $\widehat{T(q)}\simeq\nic(V)\# \ku C_{n^2}$, $V$ is
one-dimensional and $\nic_G\subseteq \nic(V)$ is generated by
$\nic_G(1)$, thus $G$ is generated by $h$ and an element $w$
subject to relations
$$ h^d=1,\;\;\; hw=q^m\; wh,\;\;\; w^n=0,$$
Where $dm=n^2$. The coaction is determined by
$$\lambda_0(h)= g^m\ot h, \quad  \lambda_0(w)=x\ot 1+ g\ot w.$$

\medbreak

Let $(A, \lambda)$ be a lifting of $G$ along $\textbf{r}(q)$. In
particular $A_0=G_0$. Using again Lemma \ref{lambda1} it is easy
to see that $A$ is generated by  elements $h$ and  $w$ with
$\lambda(h)= g^m\ot h,$ $\lambda(w)=x\ot 1+ g\ot w$ and $hw=q^m\;
wh$.

\medbreak Since $\lambda(w^n)=(1-g^n)\ot 1+ g^n\ot w^n$, whence
$w^n\in A_0$. Thus, there exists $\xi_i\in \ku$, $i=0\dots d-1$,
such that $w^n=\sum_{i=0}^{d-1} \xi_i\, h^i$. Therefore
$$\sum_{i=0}^{d-1}\xi_i\, g^n\ot h^i + (1-g^n)\ot 1= \sum_{i=0}^{d-1}\xi_i \,
g^{mi}\ot h^i.$$

This implies that $\xi_0=1$. If $m$ does not divides $n$ then
$\xi_i=0$ for all $i=1\dots d-1$, if there is an integer $a$ such
that $ma=n$ then $\xi_i=0$ for all $i\neq a$ and $\xi_a$ is
arbitrary. In the first case $A\simeq \Ac(d)$ and in the second
case $A\simeq \Bc(a,\xi_a)$.\epf

\subsection{ Module categories over $\Rep(\Hc(1,q))$} We shall assume that
$n> 2$. Let $V=\ku\{x, y\}$ denote the 2-dimensional
Yetter-Drinfeld module over $\ku C_n$, where $C_n$ is the cyclic
group generated by $g$, with action $\cdot: \ku C_n\otk V\to V$
and coaction $\delta:V\to \ku C_n\otk V$ determined by
$$ g\cdot x=q\; x,\;\;\; g\cdot y=q^{-1}\; y,\;\;\; \delta(x)=g^{-1}\ot x,\;\;\;
\delta(y)=g^{-1}\ot y.$$ It is not difficult to prove that
$\Hc(1,q)=\nic(V)\# \ku C_n$. So $\Hc(1,q)$ is a coradically
graded Hopf algebra with gradation given by
$\Hc(1,q)(i)=\nic^i(V)\otk \ku C_n$.

\begin{lema}\label{gener-in-1} Let $K=\oplus_{i=0}^{n-1} K(i)\subseteq \nic(V)$ be an
homogeneous left coideal subalgebra of $\Hc(1,q)$, that is $K$ is
graded as an algebra, $K(j) \subseteq \nic^j(V)$ and
$\Delta(K(j))\subseteq \oplus^j_{i=0} \Hc(1,q)(i)\otk K(j-i)$ for
all $j=0\dots n-1$. Then $K$ is generated as an algebra by $K(1)$.
\end{lema}

\pf If $\dim K(1)=2$ then $K= \nic(V)$ and the claim is obviously
true. Let us assume that $\dim K(1)=1$, and let $w=a x+ b y$ be a
nonzero element of $K(1)$. We will prove that for any $m=1\dots
n-1$, $K(m)$ is the 1-dimensional vector space generated by $w^m$.

Let $\theta\in K(m)$, then since $K(m)\subseteq \nic^m(V)$,
$\theta=\sum_{i=0}^m\, \alpha_i\; y^i x^{m-i} $ for some
$\alpha_i\in\ku$.

\medbreak

Let us denote by $\chi^n_i$ the quantum Gaussian coefficients (or
q-binomial coefficients), that is
$$\chi^n_i=\frac{(q^{n-i+1}-1)\dots (q^n-1)}{(q-1)\dots (q^i-1)} $$
It is well-known that
$$ \Delta(x^k)=\sum_{j=0}^k\, \chi^k_j\; g^{-j}x^{k-j} \ot x^j,
\;\;\;\Delta(y^k)=\sum_{j=0}^k\, \chi^k_j\; y^{k-j}g^{-j} \ot
y^j.$$ This implies that
\begin{align}\label{delta-w} \Delta(\theta)=\sum_{i=0}^m\sum_{j,l}\, \alpha_i\chi^{k-i}_j \chi^i_l\;
\big(y^{i-l} g^{-j-l}  x^{k-i-j}\ot y^lx^j \big).
\end{align}
Since $\Delta(\theta)\in \oplus_{i=0}^m\; \Hc(1,q)(m-i)\ot K(i)$,
there exists an element $v\in \Hc(1,q)(m-1)$ such that the summand
of $\Delta(\theta)$ that belongs to $\Hc(1,q)(m-1)\ot K(1)$ equals
$v\ot w$. Using equation \eqref{delta-w} we obtain that
$$a\,v= \sum_{i=0}^{m-1}\; \alpha_i \chi^{m-i}_1\chi^i_0\;y^i g^{-1}x^{m-i-1}
,$$
$$b\,v=\sum_{i=1}^m\; \alpha_i \chi^{m-i}_0\chi^i_1\;y^{i-1} g^{-1}x^{m-i}
.$$

When $a=0$, or $b=0$, the above equations immediately imply that
$\theta$ is a scalar multiple of $y^m$, or $x^m$ respectively. Let
us assume that $ab\neq 0$. In this case we have that
$$\frac{b}{a}\sum_{i=0}^{m-1}\; \alpha_i \;\chi^{m-i}_1\chi^i_0 \;
y^ig^{-1}x^{m-i-1}=\sum_{i=1}^m\; \alpha_i\;
\chi^{m-i}_0\chi^i_1\; y^{i-1} g^{-1}x^{m-i}.$$ Comparing
coefficients we obtain that for any $i=1\dots m-1$
$$\frac{b}{a}\,\alpha_i\chi^{m-i}_1 = \alpha_{i+1} \chi^{i+1}_1.$$
From this equation we obtain that $\alpha_i=\gamma \chi^m_i
a^{m-i}b^i$, for some $\gamma\in \ku$, thus $\theta=\gamma\,
(ax+by)^m.$\epf

I. Heckenberger pointed out to me that the above result is a very
special property that holds for quantum linear spaces, clearly not
valid for arbitrary Nichols algebras.

\medbreak

Let us define collections  of left $\Hc(1,q)$-comodule algebras
right $\Hc(1,q)$-simple. The algebras in the list are pair-wise
non isomorphic as comodule algebras. Let $d\in \Na$ be a divisor
of $n$ and $dm=n$.

\begin{itemize}
    \item[(a)] The group algebra $\ku C_d=  \ku\langle h:
    h^d=1\rangle$, with coaction determined by $\lambda(h)= g^m\ot
    h$.
    \item[(b)] For any $\xi\in \ku$, the algebras $\Ac_0(d,\xi)=
    \ku\langle h, w:  h^d=1,\, w^n=\xi 1,\, hw= q^m\; wh\rangle$ and
    coaction determined by $\lambda(h)= g^m\ot
    h$, $\lambda(w)=x\ot 1+ g^{-1}\ot w$.

 \item[(c)] For any $\xi\in \ku$, the algebras $\Ac_1(d,\xi)=
    \ku\langle h, w:  h^d=1,\, w^n=\xi 1,\, hw= q^m\; wh\rangle$ and
    coaction determined by $\lambda(h)= g^m\ot
    h$, $\lambda(w)=y\ot 1+ g^{-1}\ot w$.

    \item[(d)] For any $\xi, \mu\in \ku$, $\mu\neq 0$, the algebras
 $\Ac(\xi,\mu)=
    \ku\langle  w:   w^n=\xi 1\rangle$ with
    coaction determined by $\lambda(w)=(\mu\, x+ y)\ot 1+ g^{-1}\ot w$.

\item[(e)] If $n$ is even, then for any $\xi, \mu\in \ku$,
$\mu\neq 0$, the algebras
 $\Ac'(\xi,\mu)=
    \ku\langle  h, w: h^2=1,  w^n=\xi 1, hw=q^{\frac{n}{2}}\, wh\rangle$ with
    coaction determined by $\lambda(h)= g^{\frac{n}{2}}\ot
    h$, $\lambda(w)=(\mu\, x+ y)\ot 1+ g^{-1}\ot w$.

 \item[(f)] For any $\xi, \mu\in \ku$,  the algebras
 $\Do(d,\xi,\mu)=\ku\langle  h, w, z: h^d=1,\;  w^n=\xi 1,\; z^n=\mu
 1,\;hz=q^m\, zh,\; hw=q^{-m}\, wh,\; zw-q\, wz=0\rangle$. The coaction is determined by
 $\lambda(h)= g^m\ot h$, $\lambda(z)=x\ot 1+ g^{-1}\ot z$,
 $\lambda(w)=y\ot 1+ g^{-1}\ot w$.

\item[(g)] For any $\xi, \mu, \eta\in \ku$,  the algebras
 $\Do_1(\xi,\mu,\eta)=\ku\langle  h, w, z: h^n=1,\;  w^n=\xi 1,\; z^n=\mu
 1,\; hz=q\, zh,\; hw=q^{-1}\, wh,\; zw-q\, wz=\eta\, h^{n-2}\rangle$. The coaction is determined by
 $\lambda(h)= g\ot h$, $\lambda(z)=x\ot 1+ g^{-1}\ot z$,
 $\lambda(w)=y\ot 1+ g^{-1}\ot w$.

\item[(h)] Assume that $n=2k$, $1< k\in \Na$. For any $\xi, \mu\in
\ku$, the algebras
 $\Do_2(\xi,\mu,\eta)=\ku\langle  h, w, z: h^k=1,\;  w^n=\xi 1,\; z^n=\mu
 1,\;hz=q^2\, zh,\; hw=q^{-2}\, wh,\;  zw-q\, wz=\eta\, h^{k-1}\rangle$. The coaction is determined by
 $\lambda(h)= g^2\ot h$, $\lambda(z)=x\ot 1+ g^{-1}\ot z$,
 $\lambda(w)=y\ot 1+ g^{-1}\ot w$.

\end{itemize}

\smallbreak

\begin{lema} The algebras listed above are of the form
${}_\sigma K$, where $K\subseteq \Hc(1,q)$ is a left coideal
subalgebra and $\sigma$ is a compatible Hopf 2-cocycle. Moreover
the algebras listed above are pair-wise non Morita equivariant
equivalent.
\end{lema}
\pf The proof is completely analogous to the proof of Lemma
\ref{non-equivalent-radford}. \epf

\begin{prop}\label{modcat-for-books} Let $\Mo$ be an exact indecomposable module category
over $\Rep(\Hc(1,q))$, then $\Mo\simeq {}_K\Mo$, where $K$ is one
of the $H$-comodule algebras listed above.
\end{prop}
\pf Let $(G, \lambda_0)$ be a Loewy-graded $\Hc(1,q)$-comodule
algebra such that $G_0=\ku C_d$ for some divisor $d$ of $n$ and
$G\neq G_0$. Follows from Proposition \ref{diagram-comod2} that
$G\simeq K\# G_0$ where $K\subseteq \nic(V)$ is an homogeneous
left coideal subalgebra. Since $K(1)\neq 0$, follows from Lemma
\ref{gener-in-1} that the only four possibilities for $K$ are the
following. $K(1)=\ku<x>$,    $K(1)=\ku<y>$, $K(1)=\ku<\mu x+y>$
for some $0\neq \mu\in \ku$, or $K=\nic(V)$. Notice that if
$K(1)=\ku<\mu x+y>$, since $K$ is invariant under the action of
$G_0$ and $g\cdot x= q\; x, g\cdot y= q^{-1}\; y$ the only
possibilities are that $G_0=\ku 1$, thus $d=1$, or that $G_0\simeq
\ku C_2$.

\medbreak

Let $(A, \lambda)$ be a lifting of $G$ along $\Hc(1,q)$. In the
first three cases it is easy to prove that $A\simeq \Ac_0(d,
\xi)$, $A\simeq \Ac_1(d, \xi)$ or $A\simeq \Ac(\xi,\mu)$ for some
$\xi, \mu\in \ku$.\smallbreak

Let us assume that $(A, \lambda)$ is a lifting of $\nic(V)\# \ku
C_d$ along $\Hc(1,q)$. Using the same arguments as in the proof of
Proposition \ref{modcat-taft} one can prove that $A$ is generated
by elements $h, w,z$ subject to relations $h^d=1,\;  w^n=\xi 1,\;
z^n=\mu
 1,\;hz=q^m\, zh,\; hw=q^{-m}\, wh$ for some $\xi, \mu\in \ku$, and the
 coaction is determined by $\lambda(h)= g^m\ot h$, $\lambda(z)=x\ot 1+ g^{-1}\ot z$,
 $\lambda(w)=y\ot 1+ g^{-1}\ot w$.
Since $\lambda(zw- q wz)= g^{-2}\ot (zw- q wz)$, then $zw- q
 wz\in A_0$. Thus there exists $\eta_i\in \ku$ such that
 $zw- q wz= \sum_{i=0}^d \;\eta_i \, h^i$. Therefore
 $$\sum_{i=0}^d \;\eta_i\, g^{mi}\ot h^i=\sum_{i=0}^d
 \;\eta_i\, g^{-2}\ot h^i.$$
We conclude that there are two possibilities: if there is no
$i=0\dots d-1$ such that $g^{-2}=g^{mi}$, in which case all
$\eta_i=0$ and $zw- q wz=0$, thus $A\simeq \Do(d,\xi,\mu)$. The
other case is when there exists one $i=0\dots d-1$ such that
$g^{-2}=g^{mi}$ thus $\eta_j=0$ for all $j$ except when $j=i$.
Also $n\mid mi+2$. Since $i<d$, $mi+2< n+2$ then $n=mi+2$, but
since $m$ is a divisor of $n$, this implies that $m$ divides 2,
hence $m=1$ or $m=2$. In the first case $d=n$ and $A \simeq
\Do_1(\xi,\mu,\eta)$ and in the second case $A \simeq
\Do_2(\xi,\mu,\eta)$.\epf

\subsection{Module categories over $\Rep(u_q(\mathfrak{sl}_2))$}

For convenience we shall work over a twist equivalent Hopf algebra
rather than with $u_q(\mathfrak{sl}_2)$ itself. \smallbreak

Let $2<n\in \Na$ and $q\in \ku$ be an $n$-th primitive root of 1.
Let $\Hc_q$ be the algebra generated by elements $g, x, y$ subject
to relations
$$gx = q^2\, xg,\;  \; gy = q^{-2}\,yg,\;
 \;g^{n} = 1,\; x^n = 0,\;  \; y^n = 0,\;  \; xy-q^2\, yx =
 1-g^{-2}.$$
The algebra $\Hc_q$ is a Hopf algebra with coproduct  determined
by $\Delta(g) = g\otimes g$, $\Delta(x) = x\otimes 1 +  g^{-1}
\otimes x$, $\Delta(y) = y\otimes 1 + g^{-1} \otimes y$.

\medbreak

It is easy to verify that $\Hc_q\simeq u_q(\mathfrak{sl}_2)$ as
Hopf algebras, where the isomorphism is given by multiplying by
the group-like element $g$ and that $\gr_c \Hc_q\simeq \nic(V)\#
\ku C_n$ as Hopf algebras, where $V=\ku\{x, y\}$ is the
2-dimensional Yetter-Drinfeld module over $\ku C_n$ with structure
maps $\delta:V\to \ku C_n\otk V$, $\cdot:\ku C_n \otk V\to V$
given by
$$\delta(x)=g^{-1}\ot x,\;  \delta(y)=g^{-1}\ot y,\; g\cdot x= q^2\, x,\;
g\cdot y=q^{-2}\, y.$$

\medbreak

It is well-known that the category of representations of
$H=\nic(V)\# \ku C_n$ is tensor equivalent to $\Rep(\Hc_q)$.
Therefore we shall describe module categories over $\Rep(H)$. Let
us define a family of left $H$-comodule algebras right $H$-simple.
Let $d\in \Na$ be a divisor of $n$ and $dm=n$.

\begin{itemize}
    \item[(a)] The group algebra $\ku C_d=  \ku\langle h:
    h^d=1\rangle$, with coaction determined by $\lambda(h)= g^m\ot
    h$.
    \item[(b)] For any $\xi\in \ku$, the algebras $\Ac_0(d,\xi)=
    \ku\langle h, w:  h^d=1,\, w^n=\xi 1,\, hw= q^m\; wh\rangle$ and
    coaction determined by $\lambda(h)= g^m\ot
    h$, $\lambda(w)=x\ot 1+ g^{-1}\ot w$.
\item[(c)] For any $\xi\in \ku$, the algebras $\Ac_1(d,\xi)=
    \ku\langle h, w:  h^d=1,\,\, w^n=\xi 1,\,\, hw= q^m\; wh\rangle$ and
    coaction determined by $\lambda(h)= g^m\ot
    h$, $\lambda(w)=y\ot 1+ g^{-1}\ot w$.

     \item[(d)] For any $\xi, \mu\in \ku$, $\mu\neq 0$, the algebras
 $\Ac(\xi,\mu)=
    \ku\langle  w:   w^n=\xi 1\rangle$ with
    coaction determined by $\lambda(w)=(\mu\, x+ y)\ot 1+ g^{-1}\ot w$.

 \item[(e)] If $n$ is even, for any $\xi, \mu\in \ku$, $\mu\neq 0$, the algebras
 $\Ac'(\xi,\mu)=
    \ku\langle h, w:  h^2=1, w^n=\xi 1,  hw=q^{\frac{n}{2}}\, wh\rangle$ with
    coaction determined by $\lambda(h)=g^{\frac{n}{2}}\ot h$, $\lambda(w)=(\mu\, x+ y)\ot 1+ g^{-1}\ot w$.

 \item[(f)] If $n=4m$, for any $\xi, \mu\in \ku$, $\mu\neq 0$, the algebras
 $\Ac''(\xi,\mu)=
    \ku\langle h, w:  h^4=1,  w^n=\xi 1,  hw=q^{m}\, wh\rangle$ with
    coaction determined by $\lambda(h)=g^{m}\ot h$, $\lambda(w)=(\mu\, x+ y)\ot 1+ g^{-1}\ot w$.

    \item[(g)] For any $\xi,\mu\in \ku$, the algebras
    $\Bc(d,\xi,\mu) =\ku\langle h, w, z:  h^d=1,\,\;
    w^n=\xi 1,\,\; z^n=\mu 1,\;\, hw= q^m\; wh,\;\, hz= q^m\; zh,\;\, zw-q^2\,
    wz=1\rangle$, with coaction determined by $\lambda(h)= g^m\ot
    h$, $\lambda(w)= y\ot 1+
    g^{-1}\ot w$, $\lambda(z)= x\ot 1+ g^{-1}\ot z$.

\item[(h)] If there exists $a\in\Na$ such that $ma+2=n$, then for
any $\xi,\mu,\eta\in \ku$, with $\eta\neq 0$, set the algebras
    $\ca(d,\xi,\mu,\eta,a) =\ku\langle h, w, z:  h^d=1,\,
    w^n=\xi 1,\, z^n=\mu 1,\, hw= q^m\; wh,\, hz= q^m\; zh,\, zw-q^2\,
    wz=1+ \eta \,h^a\rangle$, with coaction determined by $\lambda(h)= g^m\ot
    h$, $\lambda(w)= y\ot 1+
    g^{-1}\ot w$, $\lambda(z)= x\ot 1+ g^{-1}\ot z$.

\end{itemize}

\begin{lema}\label{non-equi-sl2}  The algebras listed above are of the form
${}_\sigma K$, where $K\subseteq H$ is a left coideal subalgebra
and $\sigma$ is a compatible 2-cocycle. Moreover the algebras
listed above are pair-wise non Morita equivariant equivalent.
\end{lema}
\pf The proof is completely analogous to the proof of Lemma
\ref{non-equivalent-radford}. \epf

\begin{prop} If $\Mo$ is an indecomposable exact module category
over $\Rep(u_q(\mathfrak{sl}_2))$ then $\Mo$ is equivalent to one
of the following categories: ${}_{\ku C_d}\Mo$,
${}_{\Ac_0(d,\xi)}\Mo$, ${}_{\Ac_1(d,\xi)}\Mo$,
${}_{\Ac(\xi,\mu)}\Mo$, ${}_{\Ac'(\xi,\mu)}\Mo$,
${}_{\Ac''(\xi,\mu)}\Mo$, ${}_{\Bc(d,\xi,\mu)}\Mo$,
${}_{\ca(d,\xi,\mu,\eta,a)}\Mo$.
\end{prop}

\pf Let $(G, \lambda_0)$ be a Loewy-graded $\nic(V)\# \ku
C_n$-comodule algebra. We assume that $G\neq G_0$ and
$G_0=\ku\langle h: h^d=1\rangle$ for some $d$ divisor of $n$. The
space $\nic_G\subseteq \nic(V)$ is an homogeneous left coideal
subalgebra generated as an algebra by $\nic_G(1)$. This can be
proven using the same arguments as in the proof of Lemma
\ref{gener-in-1}.

Assume that $\dim \nic_G(1)=2$, then evidently $\nic_G=\nic(V)$.
Let $(A,\lambda)$ be a lifting of $(G, \lambda_0)$ along $H$.
Arguing as in Proposition \ref{modcat-for-books} we obtain that
$A$ is generated by elements $h, w,z$ subject to relations
$$h^d=1,\,\;\;
    w^n=\xi 1,\,\;\; z^n=\xi 1,\;\;\; hw= q^m\; wh,\;\;\; hz= q^m\; zh. $$
Where the coaction is determined by
 $\lambda(h)=
g^m\ot h$, $\lambda(z)=x\ot 1+ g^{-1}\ot z$, $\lambda(w)=y\ot 1+
g^{-1}\ot w$. Since
\begin{equation}\label{eq-1}\lambda(zw-q^2\, wz)= (1-g^{-2})\ot 1 +
g^{-2}\ot (zw-q^2\, wz),
\end{equation}
 then there exists $\alpha_i \in \ku$,
$i=0\dots d-1$ such that $zw-q^2\, wz= \sum_i\alpha_i\; h^i.$
Using equation \eqref{eq-1} we obtain that
$$\sum_i\alpha_i\; g^{mi}\ot h^i= (1-g^{-2})\ot 1+
\sum_i\alpha_i\; g^{-2}\ot h^i.$$ This implies that $\alpha_0=1$.
If there is no $1\leq i \leq d-1$ such that $mi=n-2$ then
$\alpha_i=0$ for all $i=1\dots d-1$ and in this case $A\simeq
\Bc(d,\xi,\mu)$. If there is a natural number $a$ such that
$ma=n-2$ then $\alpha_i=0$ for all $i\neq a$ and $\alpha_a$ is
arbitrary. In this case $A\simeq \ca(d,\xi,\mu,\alpha_a,a)$.
\medbreak

If $\dim \nic_G(1)=1$, we have three possibilities:
$\nic_G(1)=\ku< x>$, $\nic_G(1)=\ku< y>$ or $\nic_G(1)=\ku< \mu
\,x + y>$ for some $\mu\in \ku^{\times}$. In the first two cases
$A\simeq \Ac_0(d,\xi)$ or $A\simeq \Ac_1(d,\xi)$ respectively. In
the third case, since $g\cdot x= q^2\, x$ and $g\cdot y= q^{-2}\,
y$, the only possibilities are that $G_0=\ku 1$, $G_0=\ku C_2$ or
$G_0=\ku C_4$ and then $A\simeq \Ac(\xi,\mu)$, $A\simeq
\Ac'(\xi,\mu)$ or $A\simeq \Ac''(\xi,\mu)$ for appropriate $\xi,
\mu\in \ku$.

\epf

\section{Some further comments}

1. Notice that if $H$ is one of the Hopf algebras in the following
list:

\begin{itemize}
    \item[(i)] a group algebra $\ku G$ of a finite group $G$,
    \item[(ii)] a Taft Hopf algebra $T_q$,
    \item[(iii)] a Radford Hopf algebra $\textbf{r}(q)$,
    \item[(iv)] the Frobenius-Lusztig kernel
    $u_q(\mathfrak{sl}_2)$,
\end{itemize}
then exact module categories over $\Rep(H)$ are of the form
${}_{{}_{\sigma}K}\Mo$, where $K\subseteq H$ is a left coideal
subalgebra and $\sigma$ is a compatible Hopf 2-cocycle. It would
be interesting to know if this statement is true for any
finite-dimensional pointed Hopf algebra with abelian coradical.

\medbreak

2. Using the same strategy as for the book Hopf algebras and for
the Radford Hopf algebras it is likely that a classification for a
bosonization of quantum linear spaces can be worked out easily.
Also the results by Etingof and Ostrik on the classification of
exact module categories over finite supergroups $\wedge V\# \ku
G$, see \cite[section 4.2]{eo}, can be obatined using our
techniques, at least when $G$ is a cyclic group.

\medbreak

3. In principle the techniques developed here would help to
    classify exact module categories over $\Rep(u_q(\mathfrak{sl}_n))$,
    since in the paper \cite{KLS} the authors classified
    all homogeneous coideal
    subalgebras of $u_q(\mathfrak{sl}_n)$.

\subsection*{Acknowledgments}
I am grateful for the kind hospitality at LMU, specially many
thanks to Professor H.-J. Schneider. I also want to mention
Professor I. Heckenberger for sharing his insights on Nichols
algebras.

\end{document}